\documentclass[11pt,a4paper]{article}
\usepackage{fullpage}
\usepackage{amsmath,french,theorem,amssymb,epsfig}
\numberwithin{equation}{section}                

\newcommand{\cqfd}{ \nolinebreak\rule{1.5mm}{1.5mm}\par}
                    
\newcommand{\C}{{\mathbb C}}            
\newcommand{\F}{{\cal F}}
\newcommand{\id}{\mathop{\rm id}\nolimits}                      
\newcommand{\FT}{\widetilde{\cal F}}               
\newcommand{\Hom}{\mathop{\rm Hom}\nolimits}                    
\newcommand{\D}{{\cal D}}
\newcommand{\Negl}{\mathop{\rm Negl}\nolimits}                 
\newcommand{\V}{{\cal V}} 
\newcommand{\End}{\mathop{\rm End}\nolimits}                    

\newcommand{\tr}{\mathop{\rm tr}\nolimits}                     
\newcommand{\trq}{\mathop{\rm tr_q}\nolimits}                  
\newcommand{\VB}{\overline{{\cal V}}}
\renewcommand{\hom}{\mathop{\rm hom}\nolimits}                  

\renewcommand{\dim}{\mathop{\rm dim}\nolimits}                  
\newcommand{\R}{{\cal R}}
\renewcommand{\inf}{\mathop{\rm inf}\nolimits}                 
\newcommand{\RB}{\overline{{\cal R}}}
\newcommand{\CC}{{\cal C}}                                      
\newcommand{\TL}{{\cal S}}
\newcommand{\FB}{\overline{\cal F}}

\newcommand{\dessin}[1]{\vcenter{\hbox{\epsfbox{#1}}}}

\newtheorem{lem}{Lemme}[section]                             
\newtheorem{prop}[lem]{Proposition}                          
\newtheorem{theo}[lem]{Th\'eor\`eme}
\newtheorem{cor}[lem]{Corollaire}

\newcommand{\preuve}{{\bf{D\'emonstration.}} }

\newcounter{compteurenum}
\newenvironment{enum}
  {\begin{list}{\arabic{compteurenum}.}{\usecounter{compteurenum}
\setlength{\topsep}{\medskipamount}
\setlength{\partopsep}{0pt}
\setlength{\parsep}{0pt}
\setlength{\itemsep}{0pt}
\setlength{\listparindent}{\parindent}}}{\end{list}}

\newcounter{fig}
\newenvironment{fig}{\vspace*{-6mm}\begin{array}{c}\refstepcounter{fig}}
{\\ \vspace*{9mm}\text{Figure } \thefig . \end{array}\vspace*{2mm}}

\begin{document}

\title{Description topologique des repr\'esentations de $U_q(sl_2)$.}
\author{Henrik Thys \\
Institut de Recherche Math\'ematique Avanc\'ee \\
Universit\'e Louis Pasteur - CNRS \\
7, rue Ren\'e Descartes \\
67084 Strasbourg C\'edex, France \\
thys@math.u-strasbg.fr \\
}
\date{}
\maketitle

\begin{quotation}\begin{quotation} \noindent{\bf Abstract}.  We prove there
is an equivalence of ribbon categories between a category introduced by
Turaev whose morphisms are planar diagrams and a full subcategory of the
category of finite-dimensional modules over $U_q(sl_2)$.
\end{quotation}\end{quotation}

\begin{quotation}\begin{quotation} \noindent{\bf R\'esum\'e}.  Nous
construisons une \'equivalence de cat\'egories ruban\'ees entre une cat\'egorie
introduite par Turaev dont les morphismes sont des diagrammes planaires et
une sous-cat\'egorie pleine de la cat\'egorie des modules de dimension finie
sur $U_q(sl_2)$.  \end{quotation}\end{quotation}

\makeatletter\renewcommand{\@makefnmark}{}\makeatother
\footnote{\hspace*{-20pt}Mots-cl\'es :  groupe quantique, 
invariant de n\oe uds, cat\'egorie
ruban\'ee, crochet de Kauffman, alg\`ebre de Temperley-Lieb.}  
\footnotetext{\hspace*{-20pt}Classification AMS :  17B37, 57M25, 81R50.}

\section*{Introduction}

Dans cet article nous r\'esolvons un probl\`eme soulev\'e par V. Turaev
en montrant que la cat\'egorie des repr\'esentations
de $sl_2$ ou de son alg\`ebre enveloppante quantique peut se d\'ecrire
\`a l'aide de diagrammes planaires.

Le lien entre la th\'eorie des groupes quantiques et la th\'eorie des n\oe uds
est maintenant bien \'etabli. Rappelons que les groupes quantiques ont
\'et\'e introduits autour de 1983--85 par Drinfeld et Jimbo.
Ce sont des d\'eformations \`a un param\`etre des alg\`ebres enveloppantes des
alg\`ebres de Lie semisimples. La th\'eorie des repr\'esentations des groupes
quantiques peut \^etre consid\'er\'ee comme une g\'en\'eralisation puissante et
f\'econde de la th\'eorie classique des repr\'esentations des alg\`ebres de Lie
semisimples.
Les repr\'esentations des groupes quantiques forment ce que Turaev a
appel\'e
une cat\'egorie ruban\'ee, c'est-\`a-dire une cat\'egorie qui poss\`ede un produit 
tensoriel,
un tressage et une dualit\'e, ce qui permet d'\'etendre la notion de trace \`a ce 
cadre.
Ce sont ces caract\'eristiques qui expliquent qu'\`a partir d'une telle 
cat\'egorie
on peut construire des invariants d'entrelacs qui g\'en\'eralisent le 
c\'el\`ebre
invariant
construit par Vaughan Jones en 1984. La construction d'invariants topologiques
\`a partir des groupes quantiques et, plus g\'en\'eralement, de cat\'egories 
ruban\'ees
est due \`a Reshetikhin et Turaev \cite{RTu}.

Si les groupes quantiques ont des applications spectaculaires en th\'eorie
des n\oe uds
et en topologie en dimension $2$ et $3$, le probl\`eme de savoir si on peut
compl\`etement
d\'ecrire les repr\'esentations des groupes quantiques et m\^eme des groupes
classiques
en termes d'objets g\'eom\'etriques comme les n\oe uds, les tresses ou leurs
projections planaires reste largement ouvert.
Le propos de cet article est de s'attaquer au cas le plus simple,
c'est-\`a-dire \`a celui
de l'alg\`ebre de Lie semisimple $sl_2$ et de son avatar quantique.
Nous construisons ici une \'equivalence de cat\'egories entre une
cat\'egorie de
repr\'esentations de l'alg\`ebre enveloppante quantique de $sl_2$ et une
cat\'egorie $\V (a)$
construite par Turaev \cite{Tu2}, chap. XII, dans laquelle les morphismes 
sont des diagrammes planaires. Cette \'equivalence pr\'eserve les
structures ruban\'ees
qui existent aussi bien du c\^ot\'e des repr\'esentations que du c\^ot\'e 
topologique.
Nous d\'emontrons une telle \'equivalence pour toutes les valeurs 
"g\'en\'eriques"
du param\`etre $q$
qui entre dans la d\'efinition de l'alg\`ebre enveloppante quantique de $sl_2$
ainsi que lorsque $q$ est une racine de l'unit\'e d'ordre pair $\geq 6$.

Lorsque ce param\`etre est une racine de l'unit\'e diff\'erente de $1$, on sait
 que la cat\'egorie des repr\'esentations d'une alg\`ebre enveloppante quantique
 n'est plus semisimple. Lusztig a conjectur\'e une
correspondance
entre les repr\'esentations des groupes quantiques aux racines de l'unit\'e
et les repr\'esentations modulaires du groupe alg\'ebrique correspondant.
Malgr\'e la complexit\'e de la cat\'egorie des repr\'esentations pour $q$ 
racine 
de l'unit\'e, nous en obtenons n\'eanmoins une description topologique
\`a condition de nous d\'ebarrasser de ce que Turaev appelle les modules
n\'egligeables.
Ces derniers sont des modules sur lesquelles la trace est identiquement nulle
et apparaissent, par exemple, dans les travaux de Reshetikhin et Turaev 
\cite{RTu}.

La cat\'egorie $\V(a)$ de Turaev est une \'elaboration des fameuses alg\`ebres de 
Temperley-Lieb \cite{TL}
utilisant les idempotents de Jones-Wenzl \cite{Jo1,W}. Le rapport 
entre les 
alg\`ebres de Temperley-Lieb
et les invariants quantiques n'est pas nouveau. Il appara\^it par exemple
dans les travaux de \cite{GHJo,Jo2,Kf,KiM,Li1,Li2,Tu1}. 
T. Kerler \cite{kerler} a montr\'e que la cat\'egorie des repr\'esentations
de $U_q(sl_2)$ est d\'etermin\'ee par son groupe de Grothendieck, \`a ce que
Kazhdan et Wenzl \cite{KW} appelent un "twisting" pr\`es ({\em cf.} \cite{KW}
pour une extension du r\'esultat de Kerler \`a une alg\`ebre de Lie semisimple
g\'en\'erale). Comme cons\'equence, la cat\'egorie de repr\'esentations de $U_q(sl_2)$
que nous consid\'erons est \'equivalente \`a une version "tordue" de la cat\'egorie
de Turaev. Dans cet article, nous obtenons un r\'esultat plus fort qui se passe
de "twisting" en construisant explicitement une \'equivalence de cat\'egories.
 L'\'enonc\'e principal du texte est le 
th\'eor\`eme \ref{principal}.

Voici le plan de l'article. Au \S\ref{repquantsl2}, apr\`es un bref rappel
sur les cat\'egories ruban\'ees, nous r\'esumons la th\'eorie des repr\'esentations
de $U_q(sl_2)$ aussi bien dans le cas g\'en\'erique que celui des racines de 
l'unit\'e. Au
\S\ref{enonce}, nous introduisons les diagrammes planaires, la cat\'egorie
$\V(a)$ et nous \'enon\c{c}ons le th\'eor\`eme principal. 
Nous construisons un foncteur $\FT$ de $\V(a)$ vers la cat\'egorie des
repr\'esentations de $U_q(sl_2)$ au \S\ref{foncteur}.
Le \S\ref{le_4} est consacr\'e \`a une courte
\'etude des morphismes de la cat\'egorie $\V(a)$ et \`a des rappels sur la formule 
de Clebsch-Gordan. 
  La d\'emonstration que $\FT$
est une \'equivalence de cat\'egories est donn\'ee au \S\ref{preuve} pour
le cas g\'en\'erique et au \S\ref{preuvebis} pour le cas d'une racine de
l'unit\'e.  
\section{Repr\'esentations de $U_q(sl_2)$}\label{repquantsl2}
Nous  rappelons les propri\'et\'es des repr\'esentations 
de $U_q(sl_2)$ dont nous aurons besoin par la suite. Il est 
commode de les formuler en utilisant le langage des cat\'egories ruban\'ees. 
Le contenu des \S\S\ref{catmon}-\ref{mornegl} est
tir\'e de \cite{K,KRoTu,Tu2}.

\subsection{Cat\'egories mono\" idales}\label{catmon}
Rappelons qu'une {\em cat\'egorie mono\" idale} est une cat\'egorie $\CC$ munie 
d'un foncteur $\otimes : 
\CC\times \CC\rightarrow \CC$ et d'un objet
$I$, appel\'e {\em objet unit\'e}, tels que
\begin{gather}
\label{mon1}(U\otimes V)\otimes W=U\otimes (V\otimes W),  \\
\label{mon2}(f\otimes g)\otimes h=f\otimes (g\otimes h),  \\
\label{mon3}V\otimes I=I=I\otimes V,  \\
\label{mon4}f\otimes\id_I=f=\id_I\otimes f, 
\end{gather}
pour tous les objets $U$, $V$, $W$ et tous les morphismes $f$, $g$, $h$ de 
$\CC$.

Ce que nous venons de d\'efinir est usuellement appel\'e une cat\'egorie mono\" idale 
stricte. Les cat\'egories de repr\'esentations
que nous consid\'ererons ne sont \'evidemment pas strictes et il conviendrait de 
remplacer les \'egalit\'es \eqref{mon1}-\eqref{mon4}
par des isomorphismes naturels appropri\'es. L'abus que nous commettons en ne le 
faisant pas
n'est pas s\'erieux et la m\'ethode pour passer d'une cat\'egorie mono\" idale 
arbitraire \`a une cat\'egorie mono\" idale
stricte est bien connue, {\em cf.} \cite{Ml}.

Les cat\'egories mono\"idales que nous consid\'ererons par la suite sont 
$\C$-lin\'eaires, o\`u $\C$ d\'esigne
le corps des nombres complexes. Cela signifie que $\Hom_{\CC}(V,W)$ est un
 $\C$-espace 
vectoriel pour tout couple d'objets $(V,W)$ de $\CC$ et que la composition et le
 produit tensoriel des morphismes sont bilin\'eaires.

Par la suite, nous appellerons {\em cat\'egorie tensorielle} toute cat\'egorie 
mono\"idale $\C$-lin\'eaire telle que l'espace vectoriel
$\End_{\CC}(I)$ est de dimension $1$. 

Soient $({\CC},\otimes ,I)$ et $({\CC}',\otimes ,I')$ des cat\'egories 
tensorielles. Nous dirons qu'un foncteur 
${\cal G}\: :\: {\CC}\rightarrow {\CC}'$ 
{\em conserve les structures tensorielles} si les conditions suivantes sont 
r\'ealis\'ees :
\begin{enum}
\item Pour tout couple $(x,y)$ d'objets de ${\CC}$, l'application ${\cal G} 
:\Hom_{{\CC}}(x,y)\rightarrow
\Hom_{{\CC}'}({\cal G}(x),{\cal G}(y))$ est $\C$-lin\'eaire.
\item ${\cal G} (x\otimes y)={\cal G} (x)\otimes {\cal G} (y)$ pour tout couple 
$(x,y)$ d'objets et de morphismes de ${\CC}$;
\item ${\cal G} (I)=I'$ et ${\cal G}$ r\'ealise un isomorphisme de $\End_{{\CC}}
(I)$ sur $\End_{{\CC}'}(I')$;
\end{enum}

\subsection{Cat\'egories ruban\'ees}\label{catrub}
Soit $\CC$ une cat\'egorie tensorielle, d'objet unit\'e $I$. Un {\em tressage} 
est une famille 
d'isomorphismes naturels $c=(c_{V,W}:V\otimes W\rightarrow 
W\otimes V)$, o\`u $V,W$ d\'ecrivent l'ensemble des objets de $\CC$, telle que
\begin{equation*} 
c_{U\otimes V,W}=(c_{U,V}\otimes\id_V)(\id_U\otimes c_{V,W})\quad\mbox{et}\quad
c_{U,V\otimes W}=(\id_V\otimes c_{U,W})(c_{U,V}\otimes\id_W) 
\end{equation*}
pour tous les objets $U,V,W$ de $\CC$. Nous dirons que $\CC$ est munie d'une 
{\em dualit\'e} si \`a tout objet $V$ on a associ\'e
un objet $V^*$, appel\'e {\em dual} de $V$,
et des morphismes 
$b_V:I\longrightarrow V\otimes V^*$ et $d_V:V^*\otimes V\longrightarrow I$ dans
 $\CC$
tels que
\begin{equation} 
\label{axiomedualite} (id_V\otimes d_V)(b_V\otimes\id_V)=\id_V \quad\text{et}
\quad
(d_V\otimes\id_{V^*})(\id_{V^*}\otimes b_V)=\id_{V^*}.  
\end{equation}
Enfin, une {\em torsion} dans $\CC$ est une famille d'isomorphismes naturels 
$\theta =(\theta_V :V\rightarrow V)$, o\` u $V$ d\'ecrit
l'ensemble des objets de $\CC$,
telle que
\begin{equation*} \theta_{V\otimes W}=c_{V,W}c_{W,V}(\theta_V\otimes\theta_W) 
\end{equation*}
pour tous les objets $V,W$ de $\CC$.

Une {\em cat\'egorie ruban\'ee} est une cat\'egorie tensorielle $\CC$ munie d'un 
tressage $c$, d'une dualit\'e $(*,b,d)$ et 
d'une torsion $\theta$ telle que 
\begin{equation*} (\theta_V\otimes\id_{V^*})b_V=(\id_V\otimes\theta_{V^*})b_V 
\end{equation*}
pour tout objet $V$ de $\CC$.

\subsection{Morphismes n\'egligeables}\label{mornegl}
Soit $(\CC ,\otimes ,I,c,*,\theta)$ une cat\'egorie ruban\'ee, $V$ un objet de 
$\CC$ et $f$ un endomorphisme de $V$.
La {\em trace quantique} de $f$ est l'\'el\'ement de $\End_{\CC} (I)$ d\'efini 
par
\begin{equation} \tr_q (f)=d_Vc_{V,V^*}(\theta_Vf\otimes\id_{V^*})b_V. 
\end{equation}
La trace quantique jouit des propri\'et\'es suivantes :
\begin{enum}
\item Pour tous les morphismes $f:V\rightarrow W$, $g:W\rightarrow V$, nous 
avons $\trq (fg)=\trq (gf)$.
\item Pour tous les endomorphismes $f,g$ d'objets de $\CC$, nous avons $\trq
 (f\otimes g)=\trq (f)\trq (g)$.
\item Pour tout endomorphisme $k$ de $I$, nous avons $\trq (k)=k$.
\end{enum}

Un morphisme $f:V\rightarrow W$ dans $\CC$ est dit {\em n\'egligeable} si 
$\tr_q (fg)=0$ pour tout $g\in\Hom_C(W,V)$.
Un objet $V$ de ${\CC}$ est dit {\em n\'egligeable} si son morphisme identit\'e
 est n\'egligeable, {\em i.e.}
$\trq(f)=0$ pour tout endomorphisme $f$ de $V$. Il en r\'esulte que, pour tout 
objet $W$ de ${\CC}$,
tout morphisme de $\Hom_{\CC}(V,W)$ ou $\Hom_{\CC}(W,V)$ est n\'egligeable si 
$V$ est n\'egligeable. 
La composition et le produit 
tensoriel d'un morphisme n\'egligeable avec n'importe quel autre morphisme est 
un morphisme n\'egligeable. Nous en d\'eduisons que
le produit tensoriel d'un objet n\'egligeable avec n'importe quel autre objet 
est un objet n\'egligeable, {\em cf.} \cite{Tu2}, chap. XI.

Nous notons $\Negl_{\CC}(V,W)$ le sous-espace vectoriel de $\Hom_{\CC}(V,W)$ 
des morphismes n\'egligeables et 
posons $\hom_{\CC}(V,W)=\Hom_{\CC}(V,W)/
\Negl_{\CC}(V,W)$. Nous pouvons alors d\'efinir une cat\'egorie $\overline{{\CC}}$
 ayant les m\^emes objets que ${\CC}$ et telle que
$\Hom_{\overline{{\CC}}}(V,W)=\hom_{\CC}(V,W)$ pour tout couple $(V,W)$ d'objets
 de $\CC$.
Turaev  (\cite{Tu2}, chap. XI) appelle cette cat\'egorie la cat\'egorie purifi\'ee de
 ${\CC}$. Par construction, $\overline{{\CC}}$
 est une cat\'egorie ruban\'ee qui
n'admet pas de morphismes n\'egligeables non nuls.

\subsection{Rappels sur $U_q(sl_2)$}
\label{repsl2generique}

Pour plus de d\'etails, voir \cite{DeCK, CP, Ja, K, KRoTu, RTu, Tu2}. Soit 
$q$ un nombre complexe diff\'erent de $\pm 1$. L'alg\`ebre 
enveloppante quantifi\'ee $U_q(sl_2)$ est l'alg\`ebre engendr\'ee sur $\C$ par $K$, 
$K^{-1}$, $X$ et $Y$ et les relations
\begin{gather*}
KK^{-1}=K^{-1}K=1,  \\
KXK^{-1}=q^2X, \quad KYK^{-1}=q^{-2}Y,\\
XY-YX=\frac{K-K^{-1}}{q-q^{-1}}.
\end{gather*}
L'alg\`ebre $U_q(sl_2)$ est une alg\`ebre de Hopf dont la comultiplication 
$\Delta$, la co\"unit\'e $\varepsilon$
et l'antipode $S$ sont d\'etermin\'ees par
$$\begin{array}{lll}
\Delta (X)=1\otimes X+X\otimes K, & \Delta (Y)=K^{-1}\otimes Y+Y\otimes 1, & 
\Delta (K)=K\otimes K, \\
S(X)=-XK^{-1}, & S(Y)=-KY, & S(K)=K^{-1}\\
\varepsilon (X)=0, & \varepsilon (Y)=0, & \varepsilon (K)=1.
\end{array}$$ 

Nous supposerons que  $q$ est un nombre complexe g\'en\'erique (c'est \`a dire non 
racine de l'unit\'e)
 ou bien une racine primitive $2r^{\mbox{\scriptsize \`eme}}$ de
l'unit\'e, o\`u $r$ est un entier, $r\geq 3$. 
On d\'efinit $\tilde{J}=\tilde{J}
(q)$ comme l'ensemble des entiers 
$\geq 0$ si $q$ est g\'en\'erique et comme l'ensemble fini
$\{0,\dots ,r-2\}$ si $q$ est une racine primitive $2r^{\mbox{\scriptsize\`eme}}
$ de l'unit\'e.
Pour tout entier $n\in\tilde{J}$  il existe un $U_q(sl_2)$-module simple $V_n$ 
de dimension
$n+1$, isomorphe \`a son dual. 
Lorsque $q$ est g\'en\'erique, tout $U_q(sl_2)$-module simple est isomorphe \`a un 
module $V_n$, au produit tensoriel 
par un module de dimension $1$ pr\`es, {\em cf.} \cite{Lu,Ro}. Lorsque $q$ est une 
racine primitive $2r^{\mbox{\scriptsize \`eme}}$ de
l'unit\'e, tout module simple de dimension $n<r$ est isomorphe \`a $V_n$, au 
produit tensoriel par un module
de dimension $1$ pr\`es. En outre, il n'y a pas de module simple de dimension 
$n>r$ ({\em cf.} par exemple \cite{DeCK,K}). 

Notons $\R(q)$ la sous-cat\'egorie pleine des $U_q(sl_2)$-modules de dimension
finie dont les objets sont les produits tensoriels des modules simples $V_n$.
Notre but est de d\'ecrire topologiquement la cat\'egorie ruban\'ee $\RB (q)$ associ\'ee
\`a $\R(q)$ par la proc\'edure de purification d\'ecrite au \S\ref{mornegl}. 
On notera que, si $q$ est g\'en\'erique, alors $\RB(q)=\R(q)$. 
En effet, lorsque $q$ est g\'en\'erique, la cat\'egorie des $U_q(sl_2)$-modules 
de dimension finie est semisimple. Par cons\'equent, d'apr\`es \cite{Tu2}, chap. XI,
il 
n'y a pas de morphisme n\'egligeable non nul entre $U_q(sl_2)$-modules de 
dimension finie.
\section{Enonc\'e du th\'eor\`eme principal}\label{enonce}

Le but de ce paragraphe est de d\'efinir la cat\'egorie $\V(a)$ de
Turaev et d'\'enoncer le th\'eor\`eme qui la relie aux cat\'egories $\R(q)$ et 
$\RB(q)$ introduites pr\'ec\'edemment.

\subsection{Diagrammes}
\label{diagramme}
 Soient $k$ et $\ell $ deux entiers positifs. Un {\em diagramme} de type 
 $(k,\ell )$ consiste en un
nombre fini d'arcs et de cercles immerg\'es dans la bande ${\mathbb R}\times
 [0,1]$. On suppose que l'immersion n'admet
que des points doubles ordinaires et que les extr\'emit\'es des arcs sont
 $k$ points distincts fix\'es sur la droite
${\mathbb R}\times 0$ (les "entr\'ees") et $\ell $ points distincts fix\'es sur la 
droite ${\mathbb R}\times 1$ (les "sorties"). 
A chaque croisement on distingue un brin
passant dessus et un brin passant dessous. De plus, les diagrammes sont
 consid\'er\'es \`a isotopie de 
${\mathbb R}\times [0;1]$ fixant les bords pr\`es et conservant les croisements
et les passages sup\'erieurs et inf\'erieurs.

Fixons un nombre complexe $a$ non nul. On d\'efinit alors  $E_{k,\ell }$ 
comme l'espace vectoriel
sur $\C$ engendr\'e par les classes d'isotopie des diagrammes de type $(k,
\ell )$ modulo les relations suivantes :
\begin{gather}
\label{enlevercercle}{\cal D}\:\cup \:\dessin{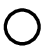}=-(a^2+a^{-2})\,
{\cal D}, \\
\label{kauf} \dessin{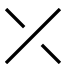}\; =a\;\dessin{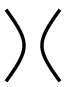}\; +a^{-1}\;
\dessin{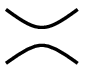}\; .
\end{gather}
Dans la premi\`ere relation, $\D\cup\dessin{cercle}$ est la classe 
d'isotopie de l'union disjointe
du diagramme $\D$ et d'un cercle. La deuxi\`eme relation lie trois 
diagrammes qui sont
\'egaux sauf dans un disque du plan o\`u ils sont comme indiqu\'e.
Cette relation s'appelle la {\em relation de Kauffman}. Il est connu 
({\em cf.} \cite{Kf}) que l'espace $E_{k,\ell }$ a pour 
base l'ensemble des classes d'isotopie des diagrammes {\em simples},
 c'est-\`a-dire ceux qui ne
contiennent ni croisements, ni cercles.

\subsection{La cat\'egorie de Temperley-Lieb}
\label{TL}

Nous sommes maintenant en mesure de d\'efinir la cat\'egorie de 
Temperley-Lieb ${\cal S}(a)$. Ses objets sont les
entiers naturels. Un morphisme $k\rightarrow \ell $ est un \'el\'ement 
de $E_{k,\ell }$. Lorsque 
${\cal D}_1$ et ${\cal D}_2$ sont des diagrammes de type respectif 
$(\ell ,m)$ et $(k,\ell )$,
la composition ${\cal D}_1\circ{\cal D}_2$ est le diagramme de type 
$(k,m)$ obtenu en
posant ${\cal D}_1$ au-dessus de ${\cal D}_2$ et en collant chaque 
sortie de $\D_2$ \`a
l'entr\'ee correspondante de $\D_1$. Cette op\'eration s'\'etend par 
lin\'earit\'e \`a tous les 
morphismes. Le morphisme identit\'e de l'objet $k$ est la classe 
d'isotopie du diagramme simple
\`a $k$ brins verticaux.

On munit $\TL(a)$ d'une structure de cat\'egorie mono\" idale de la 
mani\`ere suivante.
Le produit tensoriel de deux objets $k,\ell $ est leur somme : $k\otimes 
\ell =k+\ell $. L'objet unit\'e est $0$. Le
produit tensoriel de deux diagrammes ${\cal D}_1\otimes {\cal D}_2$ est 
donn\'e par la
juxtaposition de ces deux diagrammes, ${\cal D}_1$ \'etant plac\'e \`a 
gauche de ${\cal D}_2$. Cette
op\'eration est \'etendue \`a tous les 
morphismes par lin\'earit\'e. La cat\'egorie $\TL(a)$ est tensorielle  
car $E_{0,0}$ est
isomorphe \`a $\C$. Elle est munie
d'une structure ruban\'ee pour laquelle tout objet 
est autodual : $k^*=k$, et les morphismes de structures $c_{n,m}$, 
$\theta_n$, $b_n$ et $d_n$ sont d\'efinis par les 
diagrammes de la figure \ref{rubS(a)}.

$$\begin{fig}\label{rubS(a)}
c_{n,m}=\vcenter{\hbox{\input{tressage.pstex_t}}}, \quad \theta_n=(-1)^n\,
\vcenter{\hbox{\input{torsion.pstex_t}}},\quad
b_n=\vcenter{\hbox{\input{dualite.pstex_t}}}, \quad d_n=\vcenter{\hbox{\input{codualite.pstex_t}}}.
\end{fig}$$

\subsection{Les idempotents de Jones-Wenzl}\label{JW}
   Il est bien connu
({\em cf.} \cite{Tu2}, XII.3 et XII.4, \cite{W}) que, pour $k\geq 1$, l'alg\`ebre de 
Temperley-Lieb $E_k=E_{k,k}$
  est engendr\'ee par les $k-1$ \'el\'ements $e_1,\dots ,e_{k-1}$ o\`u 
   $e_i$ est le diagramme simple
d\'efini par la figure \ref{e_i}.  
Pour chaque $k$ strictement inf\'erieur \`a l'ordre de $a^4$,  l'alg\`ebre 
$E_k$ contient un 
idempotent central particulier, appel\'e l'idempotent de
Jones-Wenzl : il est d\'efini comme l'unique
 \'el\'ement $f_k\in E_k$ tel que $f_k-1_k$ soit un polyn\^ome non 
 commutatif sans terme constant en
  $e_1,\dots ,e_{k-1}$
 et tel que $f_ke_i=e_if_k=0$ pour tout $i=1,\dots ,k-1$.
 
$$ \begin{fig}\label{e_i}
\vcenter{\hbox{\input{e_i.pstex_t}}}
\end{fig}$$
 
\subsection{La cat\'egorie de Turaev}\label{turaev}
Nous supposons que 
le nombre complexe $a$ est g\'en\'erique ou qu'il est une racine $4r^
{\mbox{\scriptsize \`eme}}$ de 
l'unit\'e o\`u $r\geq 3$.
L'ensemble $J=J(a)$ est d\'efini comme l'ensemble des entiers $\geq 1$ 
si $a$ est g\'en\'erique et comme
l'ensemble fini $\{ 1,2,\dots ,r-2\}$ si $a$ est une racine 
$4r^{\mbox{\scriptsize \`eme}}$ de l'unit\'e.

Dans \cite{Tu2}, chap. XII, Turaev a construit une cat\'egorie
$\V(a)$ dont les objets sont les suites finies $s=(n_1,
\dots ,n_m)$ d'\'el\'ements de $J$, 
y compris la suite vide $\emptyset$. 
Pour une telle suite, on pose
$|s|=n_1+\cdots +n_m$ et
\begin{equation}
f_{s}=f_{n_1}\otimes\cdots\otimes f_{n_m}\in E_{|s|}.
\end{equation}
On convient que $f_{\emptyset}=1\in E_0$. Les morphismes de $\V(a)$ sont 
des morphismes particuliers de $\TL (a)$. 
Soient $s$ et $s'$ des objets de
$\V (a)$. Pour tout morphisme $g\in E_{|s|,|s'|}$, on pose
\begin{equation}
\widehat{g}=f_{s'}gf_{s}\in E_{|s|,|s'|}.
\end{equation}
On obtient ainsi un endomorphisme idempotent $g\mapsto\hat{g} :
E_{|s|,|s'|}\rightarrow E_{|s|,|s'|}$. Nous d\'efinissons alors
$\Hom_{\V(a)}(s ,s')$ comme l'image de cet endomorphisme :
\begin{equation}
\Hom_{\V(a)}(s ,s')= \widehat{E}_{|s|,|s'|}=\{ \widehat{g},\quad g\in 
E_{|s|,|s'|}\}.
\end{equation}
La composition des morphismes dans $\V(a)$ est induite par celle dans 
$\TL(a)$.

La cat\'egorie $\V(a)$ est une cat\'egorie ruban\'ee au sens du \S\ref{catrub}.
 Le produit tensoriel de deux
objets $s,s'$ est la concat\'enation $ss'$ des suites. L'objet unit\'e est 
la suite vide. Le produit tensoriel des morphismes
est induit par celui des morphismes de $\TL(a)$. Si $s$ et $s'$ sont des
 objets de $\V(a)$, le
tressage $c_{s,s'}:s\otimes s'\rightarrow s'\otimes s$ est donn\'e par
\begin{equation} c_{s ,s'}=\widehat{c}_{|s| ,|s'|}=(f_{s'}\otimes f_{s})
c_{|s| ,|s'|}(f_{s}\otimes f_{s'}). \end{equation}
Le dual $s^*$ d'un objet $s=(n_1,\dots ,n_m)$ est
\begin{equation} s^*=(n_m,\dots ,n_1).\end{equation}
Les morphismes $b_{s}:\C\rightarrow s\otimes s^*$, $d_{s}:s^*\otimes 
s\rightarrow\C$ et la torsion 
$\theta_{s}:s\rightarrow s$ sont d\'efinis par
\begin{gather}
 b_{s}=\widehat{b}_{|s|}=(f_{s}\otimes f_{s^*})b_{|s|}, \\ 
 d_{s}=\widehat{d}_{|s|}=d_{|s|}(f_{s^*}\otimes f_{s}), \\
 \theta_{s}=\widehat{\theta}_{|s|}=f_{s}\theta_{|s|}f_{s}.
\end{gather}

Nous noterons $\VB (a)$ la cat\'egorie ruban\'ee obtenue \`a partir de 
$\V(a)$ en annulant les morphismes n\'egligeables comme
au \S\ref{mornegl}.
Nous \'enon\c{c}ons maintenant le r\'esultat principal de l'article. C'est 
une r\'eponse au probl\`eme $24$ soulev\'e par
Turaev dans \cite{Tu2}, p. 572.
\begin{theo}
\label{principal}
Soit $a$ un nombre complexe et $q=a^2$. Alors
\begin{enum}
\item si $a$ et $q$ sont g\'en\'eriques, il existe un \'equivalence de 
cat\'egories $\FT :\V(a)\rightarrow\R(q)$;
\item si $a$ est une racine primitive $4r^{\mbox{\scriptsize \`eme}}$ 
de l'unit\'e avec $r\geq 3$ et $q$ une racine primitive 
$2r^{\mbox{\scriptsize \`eme}}$ de l'unit\'e, il existe une \'equivalence 
de cat\'egories $\FB :\VB(a)\rightarrow\RB(q)$.
\end{enum}
Les foncteurs $\FT$ et $\FB$ r\'ealisant l'\'equivalence conserve les 
structures ruban\'ees.
\end{theo}
La premi\`ere partie du th\'eor\`eme sera d\'emontr\'ee au \S\ref{preuve}, la
seconde au \S\ref{preuvebis}.
\section{Construction d'un foncteur $\V (a)\longrightarrow\R (q)$}
\label{foncteur}

Nous fixons pour toute la suite deux  nombres complexes $a$ et $q$ 
tels que $q=a^2$. Nous supposons en outre
que $a$ est soit g\'en\'erique, soit une racine primitive
 $4r^{\mbox{\scriptsize \`eme}}$ de l'unit\'e avec $r\geq 3$. 
En reprenant les notations des \S\ref{repsl2generique} et 
\S\ref{turaev}, nous avons $\tilde{J}=J\cup\{ 0\}$. 
Nous noterons $U_q$-Mod la cat\'egorie des $U_q(sl_2)$-modules de
dimension finie. 

Nous commen\c{c}ons par rappeler quelques propri\'et\'es du {\em module
fondamental} $V_1$ de dimension $2$, que nous noterons $V$.
\subsection{Le module fondamental de $U_q(sl_2)$}\label{fondamental}
Le module fondamental admet une base 
$(v_0,v_1)$ telle que
$$\begin{array}{ccc}
Kv_0=qv_0, & Xv_0=0, &Yv_0=v_1, \\
Kv_1=q^{-1}v_1, & Xv_1=v_0, &Yv_1=0.
\end{array}$$
Soit $(v^0,v^1)$ la base duale de $(v_0,v_1)$.
 Nous faisons le
choix d'une racine carr\'ee $q^{1/2}$ de $q$.
 Le tressage $c_{V,V}\in\End_{\R(q)}(V\otimes V)$ est donn\'e dans la base 
$(v_0,v_1)$ par
\begin{align}
\label{tres1}&c_{V,V}(v_0\otimes v_0)=q^{1/2}\, v_0\otimes v_0, \quad c_{V,V}
(v_0\otimes v_1)=q^{-1/2}\, v_1\otimes v_0,  \\
\label{tres2}&c_{V,V}(v_1\otimes v_1)=q^{1/2}\, v_1\otimes v_1,  \quad c_{V,V}
(v_1\otimes v_0)=q^{-1/2}\, v_0\otimes v_1+                                                                        
                                    q^{-1/2}(q-q^{-1})\, v_1\otimes v_0. 
\end{align}  
Les morphismes $b_V :\C\rightarrow V\otimes V^*$ et $d_V:V^*\otimes V
\rightarrow\C$ sont donn\'es par
\begin{equation}\label{dual}
b_V(1)=v_0\otimes v^0+v_1\otimes v^1 \quad\text{et}\quad d_V(v^i\otimes v_j)=
\delta_{ij}.
\end{equation}
Enfin la torsion $\theta_V:V\rightarrow V$ est donn\'ee par
\begin{equation}\label{torsionV}
\theta_V=q^{3/2}\,\id_V.
\end{equation}

Au \S\ref{sectionfoncteur}, nous construisons un foncteur 
$\F : \TL (a)\longrightarrow U_q$-Mod conservant les structures ruban\'ees.
La cat\'egorie $\TL(a)$ \'etant auto-duale, alors que la cat\'egorie $U_q$-Mod 
ne l'est pas (mais tout module est isomorphe \`a son
dual), nous devons faire le choix d'un isomorphisme entre le module 
fondamental $V$ et son dual.
On d\'efinit le morphisme de $U_q(sl_2)$-modules $\alpha :  V  
\rightarrow  V^* $  par 
\begin{equation}\label{alpha}
\alpha (v_0)=v^1 \quad\text{et}\quad \alpha (v_1)=-q^{-1}v^0.
\end{equation}
Il est ais\'e de voir que $\alpha$ est un isomorphisme. Posons 
\'egalement $d=d_V(\alpha\otimes\id_V):V\otimes V\rightarrow \C$ et 
$b=(\id_V\otimes\alpha^{-1})b_V:\C\rightarrow V\otimes V$, qui v\'erifient
\begin{equation}\label{b_et_d}
\begin{array}{cc}
b(1)=v_1\otimes v_0-q\, v_0\otimes v_1, & d(v_0\otimes v_0)=
d(v_1\otimes v_1)=0, \\
d(v_0\otimes v_1)=1, & d(v_1\otimes v_0)=-q^{-1}.
\end{array}\end{equation}
Des calculs \'el\'ementaires montrent alors que
\begin{gather}
\label{db} db=-(q+q^{-1}), \\
\label{tressageb_et_d} c_{V,V}=q^{1/2}\,\id_{V^{\otimes 2}}+q^{-1/2}\,bd.
\end{gather}

\subsection{Un foncteur $\F : \TL (a)\longrightarrow U_q$-Mod}
\label{sectionfoncteur}
Nous voulons d\'efinir un foncteur $\FT : \V (a)\rightarrow \R (q)$,
conservant les structures tensorielles et les structures ruban\'ees. 
Pour cela, nous rappelons
la construction bien connue d'un foncteur 
$\F : \TL (a)\longrightarrow U_q$-Mod (prop. \ref{existenceS}), 
conservant 
les structures tensorielles, ainsi que les structures ruban\'ees. 

\begin{prop}
\label{existenceS}
Il existe un unique foncteur $\F\: :\:\TL(a)\longrightarrow U_q$-Mod 
conservant les structures tensorielles tel que 
$$\F (1)=V,\quad\F \left(\dessin{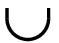}\right)=b\quad\mbox{et}\quad\F 
\left(\dessin{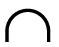}\right)=d.$$
En outre, pour tout couple $(n,m)$ d'objets de $\TL(a)$, on a
\begin{gather*}
\F(c_{n,m})=c_{\F(n),\F(m)}, \\
\F(b_n)=\left(\id_{\F(n)}\otimes (\alpha^{-1})^{\otimes n}\right)
b_{\F(n)}\quad\text{et}\quad 
   \F(d_n)=d_{\F(n)}\left(\alpha^{\otimes n}\otimes\id_{\F(n)}
   \right), \\
\F(\theta_n)=\theta_{\F(n)}.
\end{gather*}
\end{prop}
\preuve La construction d'un foncteur sur $\TL(a)$ conservant les structures tensorielles
est bien connue. Nous la rappelons bri\`evement. On d\'efinit $\F$ sur les 
objets par $\F(n)=V^{\otimes n}$. D'apr\`es
le \S\ref{diagramme}, pour d\'efinir $\F$ sur les morphismes, il suffit de
conna\^itre sa valeur sur les diagrammes simples. La valeur de
$\F(\dessin{cup})$ et $\F(\dessin{cap})$ \'etant donn\'ee, (n\'ecessairement 
$\F(\dessin{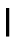})
=\id_V$ car $\F$ est un foncteur), l'image par $\F$ de tout
diagramme simple est d\'efinie par compositions et produits tensoriels. Par 
construction,
$\F$ conserve les structures tensorielles. Il reste \`a v\'erifier 
\eqref{enlevercercle}, ce qui est imm\'ediat par \eqref{db} lorsque $q=a^2$.

V\'erifions alors les quatre \'egalit\'es annonc\'ees. 
Celles concernant $\F(b_n)$ et $\F(d_n)$ sont imm\'ediates. 
Montrons que $\F(c_{n,m})=c_{\F(n),\F(m)}$, relation qu'il suffit de v\'erifier
pour $n=m=1$ gr\^ace aux propri\'et\'es des tressages.
Puisque dans $\TL(a)$ le morphisme $c_{1,1}$ est donn\'e par \eqref{kauf}, il
suffit de v\'erifier que dans $\R(q)$, on a $c_{V,V}=a\,\id_{V^{\otimes 2}}
+a^{-1}\,bd$, ce qui r\'esulte de \eqref{tressageb_et_d} 
lorsque $q=a^2$. La derni\`ere relation se montre de la m\^eme mani\`ere.\cqfd

\subsection{Le foncteur $\FT :\V (a)\longrightarrow\R (q)$}
\label{le_foncteur}
A partir du foncteur $\F :\TL(a)\rightarrow U_q$-Mod nous d\'efinissons
un foncteur $\FT : \V (a)\rightarrow U_q$-Mod comme suit :
\begin{enum}
\item $\FT ( s  )=\F (f_{ s })V^{\otimes | s |}$ pour tout objet $s$ 
de $\V(a)$;
\item $\FT (\widehat{g})=\F(\widehat{g})$ pour tout morphisme 
$\widehat{g}$ de $\V(a)$.
\end{enum}
Il est clair que $\FT$ conserve les structures tensorielles.
\begin{prop}
\label{essensurj}
Le foncteur $\FT$ est \`a valeurs dans $\R (q)$. Comme foncteur $\V (a)
\rightarrow\R (q)$, il est essentiellement surjectif.
\end{prop}
\preuve Puisque $\FT$ conserve les structures tensorielles, il suffit de 
montrer que si $n\in J$, alors $\FT(n)\cong V_n$. Soit 
$p_{n}\in\End_{\R(q)}(V^{\otimes n})$ 
l'idempotent d\'etermin\'e par 
\begin{equation*}
p_n(w)=
        \begin{cases}
        w & \text{si $w$ est de plus haut poids $n$}, \\
        0 & \text{si $w$ est de plus haut poids $p<n$}.
        \end{cases}
\end{equation*}   
Donc $p_n(V^{\otimes n})\cong V_n$. D'apr\`es \cite{M}, on a
$\F (f_n)=p_n$. Par cons\'equent, $\FT(n)=p_n(V^{\otimes n})\cong V_n$.\cqfd

Pour montrer que $\FT$ induit un foncteur $\FB:\VB(a)\rightarrow\RB(q)$
 sur les cat\'egories purifi\'ees, nous avons besoin du lemme
suivant.
\begin{lem}
\label{conservetrace}
Pour tout couple d'objets $(s,s')$ de $\V(a)$, on a
\begin{gather*}
\FT(c_{s,s'})=c_{\FT(s),\FT(s')}, \\
\FT(b_s)=\left(\id_{\FT(s)}\otimes (\alpha^{-1})^{\otimes |s|}\right)
b_{\FT(s)}\quad\text{et}\quad
\FT(d_s)=d_{\FT(s)}\left(\alpha^{\otimes |s|}\otimes\id_{\FT(s)}\right)
, \\
\FT(\theta_s)=\theta_{\FT(s)}.
\end{gather*}
Pour tout morphisme $g$ de $\V (a)$,  nous avons
$$\trq (\FT (g))=\FT (\trq (g)).$$
\end{lem}
\preuve Nous laissons le soin au lecteur de v\'erifier les quatres 
premi\`ere relations qui d\'ecoulent
de la proposition \ref{existenceS}. Nous v\'erifions la cinqui\`eme.
Soit $s$ un objet de $\V(a)$, $g\in\End_{\TL(a)}(|s|)$ et $\widehat{g}$
 l'\'el\'ement correspondant de $\End_{\V(a)}(s)$.
Par d\'efinition de la trace quantique, on a 
$$\trq (\widehat{g})=d_{ s }c_{ s  , s }(\theta_{ s }\widehat{g}\otimes
 \id_{ s })b_{ s }.$$
Par cons\'equent,
\begin{align*}
\FT (\trq (\widehat{g})) &= \FT (d_{ s })\FT (c_{ s  , s })(\FT (
\theta_{ s })\FT (\widehat{g})\otimes \FT (\id_{ s }))\FT (b_{ s }) \\
                 &= d_{\FT ( s )}(\alpha^{\otimes | s |}\otimes\id_{
                 \FT ( s  )})c_{\FT ( s ) ,\FT ( s )}
                    (\theta_{\FT ( s )}\FT (\widehat{g})\otimes \id_{
                    \FT ( s )})(\id_{\FT ( s  )}\otimes(\alpha^{-1})^
                    {\otimes | s |})b_{\FT ( s )} \\
                 &= d_{\FT ( s )}(\alpha^{\otimes | s |}\otimes\id_{
                 \FT ( s  )})c_{\FT ( s ) ,\FT ( s )}
                    (\id_{\FT ( s  )}\otimes(\alpha^{-1})^{\otimes |
                     s |})(\theta_{\FT ( s )}\FT (\widehat{g})\otimes 
                    \id_{\FT ( s )})b_{\FT ( s )} \\
                 &= d_{\FT ( s )}(\alpha^{\otimes | s |}\otimes\id_{
                 \FT ( s  )})((\alpha^{-1})^{\otimes | s |}
                 \otimes\id_{\FT ( s  )})
                    (\theta_{\FT ( s )}\FT (\widehat{g})\otimes \id_{
                    \FT ( s )})b_{\FT ( s )} \\
                 &= \trq (\FT (\widehat{g})).
\end{align*}
L'avant derni\`ere \'egalit\'e provient de la naturalit\'e du tressage. \cqfd                                 

\begin{cor}
\label{fidele}
Le foncteur $\FT$ induit un foncteur fid\`ele $\FB :\VB (a)\rightarrow
\RB(q)$.
\end{cor}
\preuve L'existence du foncteur induit est une cons\'equence imm\'ediate
 de la derni\`ere \'egalit\'e
du  lemme pr\'ec\'edent.

D\'emontrons la fid\'elit\'e de $\FB :\VB (a)\rightarrow\RB(q)$. Soit $h: s 
\rightarrow s '$ un morphisme de $\V(a)$ 
dont l'image dans $\VB(a)$ est not\'ee $\overline{h}$. Si $\FB(\overline
{h})=0$ dans $\RB(q)$,
alors $\FT(h)$ est n\'egligeable dans $\R(q)$. En particulier, pour tout 
morphisme $g:s'\rightarrow s$ de 
$\V(a)$, on a $\trq(\FT(h)\FT(g))=0$. Le lemme \ref{conservetrace} 
implique que $\FT(\trq(hg))=0$. Comme l'application
$$\FT :\End_{\V(a)}(\emptyset)\longrightarrow\End_{\R(q)}(V_0)$$
applique le diagramme vide sur l'identit\'e de $V_0$, c'est un 
isomorphisme. On en d\'eduit que $\trq(hg)=0$ ce qui prouve
que $h$ est n\'egligeable dans $\V(a)$, donc que $\overline{h}=0$.\cqfd

\section{Pr\'eliminaires techniques}\label{le_4}
Au \S\ref{sectionbontype} nous exhibons une base simple pour chaque 
espace de morphismes dans la cat\'egorie 
$\V(a)$ du \S\ref{turaev}.
 Au \S\ref{sectionCG}, nous donnons quelques rappels sur
la formule de Clebsch-Gordan pour $U_q(sl_2)$.

\subsection{Les morphismes de $\V(a)$}
\label{sectionbontype}

Etant donn\'e deux objets $ s =(n_1,\dots ,n_m)$ et $ s '=(n_1',\dots 
,n'_{m'})$ de $\V(a)$,
on dit qu'un diagramme simple $\D$  est  {\em du bon type} $( s  , s '
)$ s'il
est du type $(| s |,| s '|)$ et si
$$\widehat{\D}=f_{ s '}\D f_{ s }\not =0.$$

\begin{lem}
\label{bontype}
Les \'el\'ements $\widehat{\D}$, o\`u $\D$ est du bon type 
$( s  , s ')$, forment une base de $\Hom_{\V(a)}( s  , s ')$.
\end{lem}
\preuve Ce lemme d\'ecoule des \'egalit\'es
$f_ne_i=e_if_n=0$ pour tout $i=1,\dots ,n-1$. En effet,
pour tout $k\geq 1$, il existe un diagramme $P_k$ tel que $f_k=1_k+P_k$
et $f_kP_k=P_kf_k=0$.
Nous en d\'eduisons que
$$f_{ s }=1_{| s |}+\sum_jQ_j,\;\;\;
f_{ s '}=1_{| s '|}+\sum_{j'}Q'_{j'},$$
o\`u les $Q_j$ (resp. les $Q'_{j'}$) sont des diagrammes v\'erifiant 
$f_{ s }Q_j=Q_jf_{ s }=0$ (resp. $f_{ s '}Q'_{j'}=Q'_{j'}f_{ s '}=0$).
 Par cons\'equent, 
$$\widehat{\D}=f_{ s '}\D f_{ s }=\D+\sum_{i}\lambda_i\, \widetilde{\D}_i,$$
o\`u $\lambda_i\in\C$ et o\`u tous les $\widetilde{\D_i}$ sont des 
diagrammes simples du type $(| s |,| s '|)$ tels que $f_{ s }\widetilde{\D_i}
f_{ s '}=0$, donc tous distincts de $\D$ si $\widehat{\D}\not =0$.\cqfd

{\bf Remarque: }
Un diagramme est donc du bon type si on ne peut pas factoriser
$1^{(n_1+\cdots +n_{k-1})}\otimes e_i\otimes 
1^{(n_{k+1}+\cdots +n_m)}$ en entr\'ee ou $1^{(n'_1+\cdots +n'_{k'-1})}
\otimes e_i\otimes 
1^{(n'_{k'+1}+\cdots +n'_{m'})}$ en sortie.
La figure \ref{exemple} donne un exemple d'un diagramme du bon type 
$( s , s ')$, o\`u $ s =(3,3,3,4)$ et $ s '=(2,2)$.

$$\begin{array}{lr}
\begin{fig}\label{exemple}\dessin{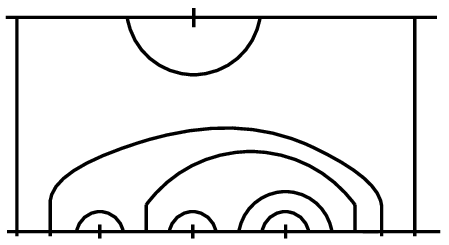}\end{fig} \qquad&
\qquad\begin{fig}\label{Dnmj}\vcenter{\hbox{\input{D_n_m_j.pstex_t}}}\end{fig}
\end{array}$$

Pour terminer ce paragraphe, nous d\'efinissons des diagrammes du bon type
particuliers. Soit $j$ un entier v\'erifiant $0\leq j\leq 
\inf (n,m)$. On d\'efinit $\D_{n,m,j}\in E_{n+m,n+m-2j}$ comme
la classe d'isotopie du diagramme simple de la figure \ref{Dnmj} o\`u,
 par convention, un entier $p$ coloriant 
un arc signifie $p$ arcs parall\`eles.

\subsection{Clebsch-Gordan}\label{sectionCG}
Le produit tensoriel des deux modules 
 $V_n$ et $V_m$ ($n,m\in \tilde{J}$) 
 v\'erifie la formule de Clebsch-Gordan
\begin{equation}
 \label{CG}
 V_n\otimes V_m\cong \!\!\!\!\!\bigoplus_{\stackrel{ k=|n-m|}{
 {\scriptscriptstyle k\equiv n+m 
 \:\mbox{\tiny mod}\, 2}}}^{n+m}\!\!\!\!\!V_k
 \end{equation}
si $n+m\in \tilde{J}$.
Si $n+m\not\in \tilde{J}$ (ce qui n'arrive que si $q$ est une racine primitive
 $2r^{\mbox{\scriptsize \`eme}}$ de l'unit\'e et $n+m>r-2$), alors
\begin{equation}
\label{CGbis}
V_n\otimes V_m\cong \!\!\!\!\!\bigoplus_{\stackrel{ k=|n-m|}
{{\scriptscriptstyle k\equiv n+m \:\mbox{\tiny mod}\, 2}}}^{2r-4-n-m}
\!\!\!\!\!V_k\oplus Z
\end{equation}
o\` u $Z$ est un module n\'egligeable. Une telle d\'ecomposition 
\'etant unique \`a isomorphisme pr\`es, nous parlerons des parties
semisimple et n\'egligeable de $V_n\otimes V_m$. ({\em cf.} \cite{And,
AndPar,CP,RTu}).

 En it\'erant \eqref{CG}, on obtient, pour tout entier $n\in J$, une 
 d\'ecomposition de la forme
 \begin{equation}
 \label{CGitere}
\displaystyle V^{\otimes n}\cong V_n\oplus \bigoplus_{
\scriptscriptstyle k_i<n}\!V_{k_i}.
\end{equation}

Nous rappelons \'egalement la forme g\'en\'erale des vecteurs de plus haut poids 
d'un produit tensoriel $W\otimes W'$
de deux modules. Soient $w\in W$ un vecteur de plus haut poids $n$, 
$w'\in W'$ un vecteur de plus haut poids $m$ et $p$
un entier $\geq 0$. Posons
\begin{equation}
\label{vecteur+hp}
v(w,w',p)=\sum_{i=0}^{p}(-1)^i\frac{[m-p+i]_q![n-i]_q!}{[i]_q![p-i]_q!
[m-p]_q![n]_q!}q^{-i(m-2p+i+1)}\, (Y^pw)\otimes (Y^{p-i}w'),
\end{equation} 
o\`u $[k]!=[k][k-1]\cdots [1]$ si $k$ est un entier $\geq 1$ et $[0]!=
1$ par convention. Nous avons alors les r\'esultats qui 
suivent ({\em cf.} \cite{K}, p.157).
\begin{enum}
\item Lorsque $q$ est g\'en\'erique, un vecteur de $W\otimes W'$ est de
 plus haut poids $k$ si et seulement s'il est \'egal \`a un vecteur $v(w,w',p)$,
  avec
$0\leq p\leq\inf (n,m)$ et $n+m-2p=k$. 
\item Lorsque $q$ est une racine primitive $2r^{\mbox{\scriptsize 
\`eme}}$ de l'unit\'e, supposons, en outre, que $W$ et $W'$ se 
d\'ecomposent en sommes directes d'une partie semisimple et n\'egligeable,
 i.e.
$$W\cong P\oplus Z,\;\; W'\cong P'\oplus Z',$$
o\`u $Z$ et $Z'$ sont des modules n\'egligeables et $P$ et $P'$ sont 
sommes directes de modules $V_n,\; n\in J$. Dans ce cas,
un vecteur de la partie semisimple de $W\otimes W'$ est 
 de plus haut poids $k$ si et seulement s'il est \'egal \`a un vecteur
 $v(w,w',p)$, avec
$0\leq p\leq\inf (n,m)$, $n+m-2p=k$ et $k\leq\inf (n+m,2r-4-n-m)$.
\end{enum}

\section{Pleine fid\'elit\'e : cas g\'en\'erique}\label{preuve}
Nous fixons un nombre complexe $a$ g\'en\'erique et posons $q=a^2$. 
Pour d\'emontrer le premier point du th\'eor\`eme \ref{principal}, il suffit,
au vu de la proposition \ref{essensurj}, d'\'etablir
la pleine fid\'elit\'e du foncteur $\FT:\V(a)\rightarrow\R(q)$  
d\'efini au \ref{le_foncteur}. Il s'agit
de montrer le th\'eor\`eme suivant.
\begin{theo}\label{theoreme}
Pour tout couple $( s , t)$ d'objets de $\V(a)$, l'application
$$\FT:\Hom_{\V(a)}( s , t)\longrightarrow\Hom_{\R(q)}(\FT( s ),\FT( 
t))$$
est un isomorphisme.
\end{theo}
Nous d\'emontrerons ce th\'eor\`eme par r\'ecurrence sur la longueur de $t$.
Les propositions \ref{casfacile} et \ref{injectivite} concernent le cas
o\`u $t$ est de longueur $1$ et le \S\ref{casgeneral} traite le cas g\'en\'eral.

Remarquons que le th\'eor\`eme est vrai lorsque $s$ et $t$ sont de 
longueur $0$ ou $1$.
En effet, les consid\'erations du \S\ref{sectionbontype} impliquent que, 
dans ce cas,
$$\Hom_{\V(a)}(s,t)=
\begin{cases}
\C\id_{s}&\text{si $s=t$}, \\
0&\text{sinon}.
\end{cases}$$
D'apr\`es le lemme de Schur, il en est de m\^eme de $\Hom_{\R(q)}(\FT(s),
\FT(t))$.

Dans toute la suite de l'article, nous conviendrons que la suite $(0)$
repr\'esente la suite vide de $\V(a)$, ce qui permet de la traiter comme
un objet de longueur $1$.
 
\subsection{Cas de $\FT :\Hom_{\V(a)}((n,m),(k))\rightarrow
\Hom_{\R(q)}(\FT ((n,m)),\FT (k))$}\label{premiercas}
On fixe deux entiers $n$ et $m$ $>0$ et un entier $k\geq 0$. Si 
$j\geq 0$ est un entier , on note $\id_j$ le morphisme identit\'e 
de $V^{\otimes j}$. Le but du \S\ref{premiercas} est de d\'emontrer le 
th\'eor\`eme \ref{theoreme} pour
$s=(n,m)$ et $t=(k)$.
\begin{lem}
\label{lem2}\label{lem4}
Si $0\leq k<|n-m|$, ou bien $n+m<k$, ou bien $k\not\equiv n+m\:\mbox
{\rm mod}\, 2$, alors
$$\Hom_{\V(a)}((n,m),(k))=0\quad \mbox{et}\quad\Hom_{\R(q)}(\FT (n,m),
\FT (k))=0.$$
Sinon,
$$\Hom_{\V(a)}((n,m),(k))=\C\,\widehat{\D}_{n,m,j},$$
(o\`u $j$ est l'unique entier v\'erifiant $0\leq j\leq \inf (n,m)$ et 
$n+m-2j=k$) et
$$\dim_{\C}\left(\Hom_{\R(q)}(\FT (n,m),\FT (k))\right)=1.$$
\end{lem}
\preuve Le calcul de la dimension de $\Hom_{\R(q)}(\FT (n,m),\FT (k))$
 provient de la formule
\eqref{CG} et du lemme de Schur. 
La d\'etermination de $\Hom_{\V(a)}((n,m),(k))$ est imm\'ediate d\`es que
l'on a remarqu\'e que s'il existe un diagramme du bon type $((n,m),(k))$,
alors c'est le diagramme $\D_{n,m,j}$, o\`u $j$ v\'erifie les conditions de
l'\'enonc\'e.\cqfd   

Nous avons maintenant besoin de la formule \eqref{vecteur+hp} et des notations 
qui y sont rattach\'ees.
Pour tout entier $j$ tel que $0\leq j\leq \inf (n,m)$, d\'efinissons
\begin{equation}\label{vnmj}
v_{n,m}^{(j)}=v(v_0^{\otimes n},v_0^{\otimes m},j)\in\FT(n,m)
\end{equation}
et
\begin{equation}\label{bnmj}
b_{n,m,j}=q^{-m+j-1}\, \frac{[n+m-j+1]_q}{[n]_q}. 
\end{equation}
D'apr\`es les rappels subordonn\'es \`a \eqref{vecteur+hp},
les vecteurs $v_{n,m}^{(j)}$ sont des vecteurs de plus haut poids 
$n+m-2j$, et pour tout
vecteur $v$ de plus haut poids $\ell $ de $\FT (n,m)$, il existe un 
(unique) $j$ avec $0\leq j\leq\inf (n,m)$, tel que $v$ soit proportionnel \`a
$v_{n,m}^{(j)}$.

\begin{lem}
\label{lem8}
Pour $n\geq 1$, $m\geq 1$ et $0\leq j\leq\inf (n,m)$, l'application
$g=\id_{n-1}\otimes d\otimes\id_{m-1}:\FT (n,m)
\rightarrow\FT (n-1,m-1)$ v\'erifie
$$g(v_{n,m}^{(j)})=b_{n,m,j}\, 
v_{n-1,m-1}^{(j-1)}$$
o\`u, par convention, $v_{n-1,m-1}^{(-1)}=0$.
\end{lem}
\preuve Rappelons que l'application $d:V\otimes V\rightarrow\C$ est d\'efinie 
par \eqref{b_et_d}. On remarque que 
$g(v_{n,m}^{(j)})\in\FT (n-1,m-1)$ et est de plus haut poids $n+m-2j=
(n-1)+(m-1)-2(j-1)$. On peut alors utiliser la remarque pr\'ec\'edant 
l'\'enonc\'e: on a
$g(v_{n,m}^{(j)})=\gamma\, v_{n-1,m-1}^{(j-1)}$,
o\`u $\gamma$ est un certain nombre complexe. Pour calculer $\gamma$, il suffit
par exemple de comparer les coefficients devant $v_0^
{\otimes (n-1)}\otimes Y^{j-1}v_0^{\otimes (m-1)}$ dans $g(v_{n,m}^{(j)})$ et
$v_{n-1,m-1}^{(j-1)}$. Comme 
$$\Delta (Y^i)=\sum_{s=0}^{i}q^{s(i-s)}
\frac{[i]_q!}{[s]_q![i-s]_q!}(Y^sK^{s-i})\otimes Y^{i-s}$$
 ({\em cf.} par
exemple \cite{K}), nous avons 
$$g\left(Y^iv_0^{\otimes n}\otimes Y^jv_0^{\otimes m}\right)=
\delta_{i0}[j]_qY^iv_0^{\otimes (n-1)}\otimes Y^{j-1}v_0^{\otimes (m-1)}
-\delta_{j0}[i]_qq^{i-j-n-1}Y^{i-1}v_0^{\otimes (n-1)}\otimes Y^{j}v_0^
{\otimes (m-1)},$$
o\`u $\delta_{kl}$ est le symbole de Kronecker.
Posons alors
$$c_{j,i}^{n,m}=(-1)^i\frac{[m-j+i]_q![n-i]_q!}{[i]_q![j-i]_q![m-j]_q!
[n]_q!}q^{-i(m-2j+i+1)}.$$
Nous avons donc 
$$g(v_{n,m}^{(j)})=(c_{j,0}^{n,m}[j]_q-c_{j,1}^{n,m}q^{-j-n+1})\; v_0^{\otimes 
(n-1)}\otimes Y^{j-1}v_0^{\otimes (m-1)}+v,$$
o\`u $v$ est une combinaison lin\'eaire de $Y^iv_0^{\otimes (n-1)}\otimes 
Y^{j-1-i}v_0^{\otimes (m-1)}$ avec $i\not =0$.
Or
\begin{align*}
c_{j,0}^{n,m}[j]_q-c_{j,1}^{n,m}q^{-j-n+1}&=\frac{[m-j]_q![n]_q![j]_q}
{[j]_q![m-j]_q![n]_q!}+\frac{[m-j+1]_q![n-1]_q!}
                                                    {[j-1]_q![m-j]_q!
                                                    [n]_q!}q^{-m+2j-2-n-j+1} \\
&=\frac{[n]_q+[m-j+1]_qq^{-m-n+j-1}}{[j-1]_q![n]_q}=b_{n,m,j}\, c_{j-1,0}^{n-1,m-1},
\end{align*}
ce qui ach\`eve la d\'emonstration du lemme.\cqfd

\begin{prop}
\label{casfacile}
L'application $\FT :\Hom_{\V(a)}((n,m),(k))\longrightarrow\Hom_{\R(q)}
(\FT((n,m),(k)))$ est un isomorphisme.
\end{prop}
\preuve Au vu du lemme \ref{lem2}, il nous reste \`a traiter le cas 
$$|n-m|\leq k\leq n+m\:\mbox{ et }\: k\equiv n+m\:\mbox{mod}\, 2,$$
{\em i.e.} le cas o\`u il existe un entier $j$ tel que $0\leq j\leq
\inf(n,m)$ et $n+m-2j=k$. D'apr\`es ce m\^eme lemme, les deux espaces sont
de dimension $1$; donc il suffit de montrer que $\FT(\D_{n,m,j})\not =0$.
Comme 
$$\F(\D_{n,m,j})=(\id_{n-j}\otimes d\otimes\id_{m-j})\cdots (\id_{n-1}
\otimes d\otimes\id_{m-1})$$
et $\F(\D_{n,m,0})=\id_{n+m}$, on obtient, 
en appliquant le lemme \ref{lem8},
$$\F (\D_{n,m,j})(v_{n,m}^{(j)})=b_{n,m,j}\cdots b_{n-j+1,m-j+1,1}\, 
v_{n-j,m-j}^{(0)}.$$
Or, pour $\widehat{\D}_{n,m,j}=f_k\,\D_{n,m,j}(f_n\otimes f_m)$, on a
\begin{gather*}
\FT (\widehat{\D}_{n,m,j})=f_k\F (\D_{n,m,j})(f_n\otimes f_m), \\
(f_n\otimes f_m)(v_{n,m}^{(j)})=v_{n,m}^{(j)}, \\
f_k(v_0^{\otimes k})=v_0^{\otimes k}, \\
v_{n-j,m-j}^{(0)}=v_0^{\otimes k},
\end{gather*}
car $v_{n,m}^{(j)}\in \FT (n,m)$ et $v_0^{\otimes k}\in \FT (k)$. Comme
 tous les $b_{n,m,j}$ sont non nuls, 
$\FT(\D_{n,m,j})(v_{n,m}^{(j)})\not =0$.\cqfd

\subsection{Cas de $\mathbf \FT :\Hom_{\V(a)}( s ,(k)) \longrightarrow
\Hom_{\R(q)}(\FT ( s ),\FT (k))$}\label{deuxiemecas}

Nous allons maintenant \'etablir la pleine fid\'elit\'e dans le cas o\`u
 $ s $ est une suite de longueur quelconque.
 Fixons un entier $m\geq 3$ et des entiers strictement positifs
$n_1,\dots ,n_m$. Nous notons $ s $ la suite $(n_1,\dots ,n_m)$ et 
$ s '$  la suite tronqu\'ee $(n_1,\dots ,n_{m-1})$. 
Nous fixons \'egalement un entier positif $k$. Nous appellerons {\em espace 
de plus haut poids} $j$ d'un module $W$ 
le sous-espace vectoriel engendr\'e par les vecteurs de $W$ de plus haut 
poids $j$.

Soit $\{ b_1,\dots ,b_N\}$ l'ensemble des entiers tel que
pour tout $j=1,\dots ,N$, la dimension $p_j$ de l'espace de plus haut 
poids $b_j$ de $\FT(s')$ soit strictement positive et tel qu'il existe un 
entier $\ell _j$ v\'erifiant $0\leq \ell _j\leq \inf (b_j,n_m)$ et 
$b_j+n_m-2\ell _j=k$. (Les notations ne signifient pas que cet ensemble 
est non vide). Les deux lemmes suivants justifient la donn\'ee de 
cet ensemble.
 
\begin{lem}
\label{lem11}
Supposons qu'il existe un vecteur de plus haut poids $k$ dans  $\FT 
( s )$.
Alors l'ensemble $\{ b_1,\dots ,b_N\}$ est non vide.
\end{lem}
\preuve Soit $v$ un vecteur de plus haut poids $k$ dans $\FT ( s )$.
Nous avons
$$\FT ( s )\cong \FT ( s ')\otimes \FT (n_m).$$
D'apr\`es la formule \eqref{vecteur+hp}, il existe un vecteur $w$ de
 plus haut poids $b$ et un entier $\ell $ tel que
$v=v(w,v_0^{\otimes n_m},\ell )$,
avec $0\leq \ell \leq\inf (b,n_m)$ et $b+n_m-2\ell =k$.\cqfd

Pour tout $j=1,\dots ,N$, soit $w_1^j,\dots ,w_{p_j}^j$ une base de l'espace
de plus haut poids $b_j$ de $\FT(s')$.

\begin{lem}
\label{lem12}
Supposons qu'il existe un vecteur de plus haut poids $k$ dans $\FT (s)$. 
Alors la famille $(v_p^{j})=v(w_p^{j},v_0^{\otimes n_m},\ell _j)$,
$j=1,\dots ,N$ et $p=1,\dots ,p_j$, forme 
une base de l'espace de plus haut poids $k$ de $\FT(s)$.
\end{lem}
\preuve Remarquons tout d'abord que ces vecteurs sont bien dans $\FT 
( s  )=\FT ( s ')\otimes \FT (n_m)$ puisque $w_p^{j}\in \FT ( s ')$
et $v_0^{\otimes n_m}\in\FT (n_m)$. Avec les notations rattach\'ees \`a la
 formule \eqref{vecteur+hp}, posons $W=\FT ( s ')$
et $W'=\FT (n_m)$. Puisque $0\leq j\leq\inf(b_j,n_m)$ et $k=b_j+n_m-2
\ell _j$, les vecteurs $v_p^{j}$ sont 
des vecteurs de plus haut poids $k$ de $\FT ( s  )$. Ils sont 
lin\'eairement ind\'ependants car les vecteurs
$(Y^iw_p^{j})$ (o\`u $i=1,\dots ,\ell _j$, $p=1,\dots ,p_j$ et $j=1,
\dots ,N$) le sont
par construction. 

Soit alors $v\in\FT ( s  )$ un vecteur de plus haut poids $k$. Puisque
 $W'\cong V_{n_m}$, tous ses vecteurs de plus haut poids
sont proportionnels \`a $v_0^{\otimes n_m}$ et de poids $n_m$, et donc, 
d'apr\`es la formule \eqref{vecteur+hp}, il existe $w\in W$
de plus haut poids $b$ tel que $v=v(w,v_0^{\otimes n_m},\ell )$ et tel
 que $0\leq \ell \leq\inf (b,n_m)$ et $k=b+n_m+\ell $.
Donc $b=b_j$ et $\ell =\ell _j$ pour un certain $j$; donc $w=\sum_{p=1}
^{p_j}\lambda_p\, w_p^{j}$, $\lambda_p\in\C$,
et $v=\sum_{p=1}^{p_j}\lambda_p\, v_p^{j}$, ce qui prouve que les 
vecteurs $v_p^{j}$ forment une base de l'espace de plus haut poids $k$
de $\FT(s)$.\cqfd

Pour tout $j=1,\dots ,N$, on note $\{\D_{j,p}\}_{p=1,\dots ,p'_j}$ 
l'ensemble des diagrammes du bon type $( s ',(b_j))$. On pose
\begin{equation}\label{DB}
\overline{\D_{j,p}}=\D_{b_j,n_m,l_j}(\D_{j,p}\otimes 1_{n_m}),
\end{equation}
$j=1,\dots N$, $p=1,\dots ,p'_j$, o\`u $\D_{b_j,n_m,l_j}$ est d\'efini 
par la figure \ref{Dnmj}. 

\begin{prop}
\label{recdiag}
Soit $k\geq 0$. Alors l'ensemble des diagrammes $\overline{\D_{j,p}},\;j=1,
\dots N,\; p=1,\dots ,p'_j$ forme une base de $\Hom_{\V(a)}(s,(k))$.
\end{prop}
\preuve Nous allons montrer que les diagrammes $\overline{\D_{j,p}}$ sont 
distincts et 
du bon type $(s,(k))$ et que tout diagramme du bon type $(s,(k))$ est \'egal 
\`a un diagramme $\overline{\D_{j,p}}$. On conclura par le lemme \ref{bontype}.

Soit $\D$ un diagramme du bon type $(s,(k))$. Nous allons 
le transformer de la mani\`ere suivante : on remonte tous les brins 
s'attachant au bloc $n_m$ des entr\'ees (et uniquement
ceux-ci), puis on divise $\D$ en deux diagrammes $\D=\D_1\circ \D_2$,
qui sont uniquement d\'etermin\'es par $\D$.
La figure \ref{construction} donne l'id\'ee de la construction.

$$\begin{fig}\label{construction}
\vcenter{\hbox{\input{recdiag1.pstex_t}}}=\vcenter{\hbox{\input{recdiag2.pstex_t}}}
\end{fig}$$

Celle-ci implique tout d'abord que $\D_2=\D'\otimes 1_{n_m}$.
Le diagramme $\D'$ est du bon type $(s',(b))$ o\`u $b\geq 0$ est un 
entier non encore d\'etermin\'e. En effet, les
entr\'ees de $\D'$ sont celles de $\D$. En ce qui concerne les sorties,
 il faut montrer
qu'il n'y a pas d'arc les reliant. Or les brins s'attachant aux sorties
 de $\D'$ proviennent soit d'un brin reliant une
entr\'ee et une sortie de $\D$, soit d'un arc reliant une entr\'ee de 
$\D$ \`a une entr\'ee du bloc $n_m$, arc que l'on a remont\'e
par la construction pr\'ec\'edente. Dans les deux cas, ces brins ne 
peuvent relier deux sorties de $\D'$.
De la m\^eme mani\`ere, on montre que $\D_1$ est du bon type $((b,n_m),
(k))$. D'apr\`es le lemme \ref{lem4}, 
on sait alors que $\D_1=\D_{b,n_m,\ell }$
o\`u $\ell $ est un entier v\'erifiant $b+n_m-2\ell =k$ et $0\leq \ell
 \leq\inf (b,n_m)$. Par d\'efinition des $b_j$, il existe $j$ tel
que $b=b_j$ et $\ell =\ell _j$. 

Il reste \`a prouver que $\overline{\D_{j,p}}$ et $\overline{\D_{j',p'}}$
 sont distincts d\`es que $j\not =j'$ ou
$p\not =p'$. Notons $t$ (resp. $t'$) le nombre de brins du bloc 
d'entr\'ees $n_m$ de $\overline{\D_{j,p}}$ (resp. $\overline{\D_{j',p'}}$) 
reliant des sorties. Alors $t=n_m-j$ et $t'=n_m-j'$; donc si $j\not =j'
$, alors $\overline{\D_{j,p}}\not =\overline{\D_{j',p'}}$.
Si $j=j'$, alors $p\not =p'$ et donc $\D_{j,p}\not =\D_{j,p'}$, ce qui
prouve la proposition. \cqfd

Nous allons maintenant d\'emontrer deux propositions qui prouvent le
th\'eor\`eme \ref{theoreme} dans le cas o\`u $t$ est de longueur $1$. 
Rappelons pour cela que tout morphisme 
de $\Hom_{\R(q)}(\FT (s),\FT (k))$ est uniquement d\'etermin\'e par ses 
valeurs sur les vecteurs 
de plus haut poids de $\FT (s)$, et qu'il envoie tout vecteur de plus 
haut poids $k'\not =k$ sur $0$. En particulier, la dimension de 
$\Hom_{\R(q)}(\FT (s),\FT (k))$ est \'egale \`a la dimension de l'espace
de plus haut poids $k$ de $\FT (s)$.

\begin{prop}
\label{egalitedim}
Pour tout entier $k\geq 0$, on a
$$\dim \left(\Hom_{\V (a)}(s,(k))\right)=\dim\left(\Hom_{\R(q)}(\FT (s),
\FT (k))\right).$$
\end{prop}
\preuve D'apr\`es la remarque pr\'ec\'edente et le lemme \ref{bontype}, il 
suffit de 
d\'emontrer que le nombre de diagramme du bon type $(s,(k))$
est \'egal \`a la dimension de l'espace de plus haut poids
$k$ de $\FT (s)$. Nous allons d\'emontrer ceci par r\'ecurrence sur la 
longueur $m$ de $s$. Le cas $m=2$ est l'objet du lemme \ref{lem4}.
Supposons alors que pour tout $\widetilde{m}<m$, tout entier 
$\widetilde{k}\geq 0$ et tout objet $\widetilde{ s }$ 
de longueur $\widetilde{m}$, le nombre de diagramme
du bon type $(\widetilde{ s },(\widetilde{k}))$ est \'egal \`a la 
dimension de l'espace de plus haut poids 
$\widetilde{k}$ de $\FT (\widetilde{ s })$.

Montrons d'abord que l'existence d'un diagramme $\D$ du bon type 
$( s ,(k))$ implique l'existence d'un vecteur de plus haut poids $k$
dans $\FT ( s )$. Nous savons d'apr\`es la proposition \ref{recdiag} que 
$\D=\overline{\D_{j,p}}$ pour un unique
couple $(j,p)$, et  donc il existe un diagramme $\D_{j,p}$ du bon 
type $(s',(b_j))$. Puisque
$s'$ est de longueur $m-1$, nous savons qu'il existe un vecteur 
$w\in\FT ( s ')$
de plus haut poids $b_j$ par hypoth\`ese de r\'ecurrence. Alors le 
vecteur $v(w,v_0^{\otimes n_m},\ell _j)$ est de plus 
haut poids $b_j+n_m-2\ell _j=k$ d'apr\`es \eqref{vecteur+hp}.

R\'eciproquement, si  $\FT ( s )$ admet un vecteur $v$ de plus haut poids
$k$, il existe, d'apr\`es la formule \eqref{vecteur+hp}, un vecteur $w\in\FT(s')$
de plus haut poids $b$ et un entier $\ell$ tel que 
$b+n_m-2\ell =k$, tels que $v=v(w,v_0^{\otimes n_m},\ell )$.
Par hypoth\`ese de r\'ecurrence, il 
existe un diagramme $\D'$ du bon type $(s',(b))$. Puisque le diagramme
 $\D_{b,n_m,\ell }$ est bien d\'efini (car $b+n_m-2\ell =k$), le
diagramme $\D_{b,n_m,\ell }(\D'\otimes 1_{n_m})$ \'egalement et est du 
bon type $(s,(k))$.

Nous venons de montrer qu'il n'y a pas de diagramme du bon type $(s,(k))$ 
si et seulement si la dimension de l'espace de plus haut poids $k$ de 
$\FT(s)$ est nulle. Supposons maintenant cette dimension $>0$. Reprenons 
les notations que nous avons utilis\'ees pour
d\'efinir les diagrammes $\overline{\D_{j,p}}$. D'apr\`es la proposition
\ref{recdiag} le nombre de diagrammes du bon type $( s  ,(k))$ est 
$p'_1+\cdots +p'_q$. D'apr\`es le lemme \ref{lem12}, la dimension de 
l'espace de plus haut poids $k$ de $\FT ( s  )$ est $p_1+\cdots +p_q$. Or, par 
hypoth\`ese de r\'ecurrence, nous avons $p'_j=p_j$, 
$1\leq j\leq q$ puisque $p_j$ est la dimension de l'espace de plus
 haut poids $b_j$ de $\FT ( s ')$, ce qui prouve le lemme.\cqfd

Le th\'eor\`eme \ref{theoreme} est donc d\'emontr\'e lorsque $t$ est de
longueur $1$ et si $\FT (s)$ ne contient
pas de vecteur de plus haut poids $k$. Nous supposerons donc jusqu'\`a la
fin du \S\ref{preuve} que $\FT (s)$ contient un vecteur de plus haut poids
$k$.

\begin{prop}
\label{prop16}\label{injectivite}
L'application 
$\FT :\Hom_{\V (a)}(s,(k))\rightarrow\Hom_{\R(q)}(\FT(s),\FT(k))$
est un isomorphisme.
\end{prop}
\preuve Notons $g$ la dimension commune de $\Hom_{\V (a)}(s,(k))$ et
$\Hom_{\R(q)}(\FT(s),\FT(k))$. D'apr\`es la remarque pr\'ec\'edant la
proposition \ref{egalitedim}, il suffit de prouver qu'il existe une base
$v_1,\dots ,v_g$ de l'espace de plus haut poids $k$ de $\FT(s)$ telle que,
si $\D_1,\dots ,\D_g$ sont les diagrammes du bon type $(s,(k))$, il 
existe des nombres complexes $\alpha_1,\dots ,\alpha_g,$ 
tous non nuls tels que 
\begin{align}
\F (\D_1)(v_1)&=\alpha_1\, v_0^{\otimes k}, \\
\F (\D_p)(v_{p'})&=0\quad\text{si $1\leq p'<p\leq g$}, \\
\F (\D_p)(v_p)&=\alpha_p\, v_0^{\otimes k}\quad\text{si $p\leq g$}, 
\end{align}
o\`u $v_0^{\otimes k}\in V^{\otimes k}$.
Nous d\'emontrons ceci par r\'ecurrence sur la longueur $m$ de $s$. 
Le cas $m=2$ est l'objet 
du \S\ref{premiercas} o\`u nous avons montr\'e que $q=1$, $p_1=1$, $g=1$ 
et que
$$\F(D_{n,m,j})(v_{n,m}^{(j)})=b_{n,m,j}\cdots b_{n-j+1,m-j+1,1}
v_0^{k}.$$ 

Supposons $m\geq 3$. Pour tout $j=1,\dots N$, nous choisissons une base
$w_p^{j}$,  $1\leq p\leq p_j$, de l'espace de plus haut poids $b_j$ de 
$\FT(s')$ telle que,
pour tout $j=1,\dots ,N$, il existe, par hypoth\`ese de r\'ecurrence, des
 nombres complexes 
$\alpha_{j,p}$ tous non nuls ($1\leq p\leq p_j$) tels que
\begin{align}
\F (\D_{j,1})(w_1^{j})&=\alpha_{j,1}\, v_0^{\otimes b_j}, \\
\label{hyprecprop16a}\F (\D_{j,p})(w_{p'}^{j})&=0\quad\text{si 
$1\leq p'<p\leq p_j$}, \\
\label{hyprecprop16b}\F (\D_{j,p})(w_p^{j})&=\alpha_{j,p}\, 
v_0^{\otimes b_j},\quad\text{si $p\leq p_j$}. 
\end{align}
Fixons alors un $j$ et soient $p$ et $p'$ des entiers tel que 
$1\leq p,p'\leq p_j$. On a
\begin{align*}
\F (\overline{\D_{j,p}})(v_{p'}^{j})&=\F (\D_{b_j,n_m,\ell _j})
\F (\D_{j,p}\otimes 1_{n_m})(v_{p'}^{j}) \\
&=\F (\D_{b_j,n_m,\ell _j})\left(\sum_{i=0}^{\ell _j}c_{\ell _j,i}^
{b_j,n_m}\, (Y^i\F (\D_{j,p})(w_{p'}^{j}))\otimes 
                                                (Y^{\ell _j-i}v_0^{
                                                \otimes n_m})\right) 
\end{align*}
o\`u la deuxi\`eme \'egalit\'e provient de la formule \eqref{vecteur+hp} et
 de la d\'efinition du vecteur $v_{p'}^{j}$
({\em cf.} lemme \ref{lem12}). Si $p'<p$, la formule \eqref{hyprecprop16a}
implique que  $\F (\overline{\D_{j,p}})(v_{p'}^{j})=0$. Si $p=p'$, 
on a, d'apr\`es \eqref{vecteur+hp} et \eqref{hyprecprop16b}, 
$$\F (\overline{\D_{j,p}})(v_{p'}^{j})=\F (\D_{b_j,n_m,\ell _j})
\left(\alpha_{j,p}\sum_{i=0}^{\ell _j}c_{\ell _j,i}^{b_j,n_m}\, 
         (Y^iv_0^{\otimes b_j})\otimes (Y^{\ell _j-i}v_0^{\otimes n_m})\right)
         =v(v_0^{\otimes b_j}, v_0^{\otimes n_m} ,\ell _j),$$
qui est donc 
un vecteur de plus haut poids $k$ de $\FT (s)$. Or, d'apr\`es 
le \S\ref{premiercas} et le lemme \ref{lem4}, il existe
$\gamma_j\not =0$ tel que 
$$\F (\D_{b_j,n_m,\ell _j})(v(v_0^{\otimes b_j},v_0^{\otimes n_m},
\ell _j))=\gamma_jv_0^{\otimes k}.$$
On en d\'eduit que
$$\F (\overline{\D_{j,p}})(v_{p'}^{j})=
 \begin{cases}
   0&\text{si $p'<p$}, \\
          \alpha_{j,p}\gamma_j\, v_0^{\otimes k}&\text{si $p'=p$}.
         \end{cases}
  $$
De plus, on peut toujours supposer $b_1>\cdots >b_N$. Par cons\'equent, 
$\F (\D_{j,p})(w_{p'}^{j'})=0$,
pour tous $j'<j$, $1\leq p\leq p_j$ et $1\leq p'\leq p_{j'}$.
En effet, rappelons tout d'abord que $w_{p'}^{j'}\in\FT ( s ')$ ($1
\leq j'\leq N,\:1\leq p'\leq p_{j'}$)
et que $\D_{j,p}\in E_{| s '|,b_j}$ ($1\leq j\leq N,\:1\leq p\leq 
p_{j}$). Donc 
$\F (\D_{j,p})(w_{p'}^{j'})$ a bien un sens. Ensuite, remarquons que
 $w_{p'}^{j'}$ est de plus haut poids $b_{j'}>b_j$. 
Donc $\F (\D_{j,p})(w_{p'}^{j'})\in V^{\otimes b_j}$
l'est \'egalement. Mais rappelons ({\em cf.} \eqref{CGitere}) que tous 
les vecteurs de plus haut poids de $V^{\otimes b_j}$ ont
un poids inf\'erieur \`a $b_j$, ce qui prouve que $\F (\D_{j,p})(w_{p'}^
{j'})=0$.

En renum\'erotant alors les vecteurs $\{ v_p^{j}\}_{\stackrel{p=1,
\dots p_j}{\scriptscriptstyle j=1,\dots ,N}}$
et les diagrammes $\{\overline{\D_{j,p}}\}_{\stackrel{p=1,\dots p_j}{
\scriptscriptstyle j=1,\dots ,N}}$
selon l'ordre
$$v_1^{1},\: v_2^{1},\dots ,v_{p_1}^{1},\: v_1^{2},\dots ,v_{
p_2}^{2},\dots ,\dots ,v_p^{j},\: v_{p+1}^{j},
\dots ,v_1^{N},\dots ,v_{p_N}^{N},$$
(et l'ordre correspondant pour les diagrammes), nous obtenons deux 
familles satisfaisant aux 
conditions de l'\'enonc\'e.\cqfd

\subsection{D\'emonstration du th\'eor\`eme \ref{theoreme}}\label{casgeneral}
D'apr\`es \cite{K}, XIV.2.2, pour tout triplet d'objets $(U,V,W)$ d'une 
cat\'egorie ruban\'ee, l'application
$\sharp :\Hom (U\otimes V,W)\longrightarrow\Hom (U,W\otimes V^*) $
donn\'ee par
$$f^{\sharp}=(f\otimes\id_V)(\id_U\otimes b_V),$$
et l'application
$\flat :\Hom (U,W\otimes V^*)\longrightarrow\Hom (U\otimes V,W^*)$
donn\'ee par
$$g^{\flat}=(\id_W\otimes b_V)(g\otimes\id_V),$$
sont des bijections inverses l'une de l'autre. Dans le cas des 
cat\'egories $\V(a)$ et $\R(q)$, ce sont 
en outre des 
isomorphismes d'espaces vectoriels. Soit
alors $s$ et $t$ deux objets quelconques de $\V(a)$. On suppose que 
$t$ est de longueur $>1$.
Soit $k>0$ l'entier et $\tilde{s}$ la suite tels que $t=\tilde{s}\otimes (k)$.
Nous allons d\'emontrer le th\'eor\`eme \ref{theoreme} par r\'ecurrence sur 
la longueur de la suite $t$. 
Le cas o\`u $t$ est de longueur $1$ (ou $0$) a \'et\'e d\'emontr\'e au 
\S\ref{deuxiemecas}. 
D'apr\`es l'hypoth\`ese de r\'ecurrence, l'application
$$\FT :\Hom_{\V(a)}(s\otimes (k), \tilde{s} )\longrightarrow\Hom_{\R(q)}
\left(\FT (s\otimes k),\FT (\tilde{s})\right)$$
est un isomorphisme. Par cons\'equent, l'application qui \`a $\widehat{\D}
\in\Hom_{\V(a)}(s,t)$ associe
$$\left(\FT ((\widehat{\D})^{\flat})\right)^{\sharp}\in\Hom_{\R(q)}
\left(\FT (s),\FT (\widetilde{s})\otimes\FT (k)^*\right)$$
est un isomorphisme. Or 
\begin{align*}
\FT ((\widehat{\D})^{\flat})&=\left(\id_{\FT ( \widetilde{s} )}\otimes
\FT (\widehat{d}_k)\right)
                                          \left(\FT (\widehat{\D})
                                     \otimes\id_{\FT (k)}\right) \\
&=\left(\id_{\FT ( \widetilde{s} )}\otimes (d_{\FT (k)}(\alpha^{
\otimes k}\otimes\id_{\FT (k)}))\right)\left(
           (\FT (\widehat{\D})\otimes\id_{\FT (k)})\right),
\end{align*}    
o\`u la deuxi\`eme \'egalit\'e provient du lemme \ref{conservetrace}.                                                                        
En utilisant la repr\'esentation graphique des morphismes dans une 
cat\'egorie ruban\'ee (donn\'ee, par
exemple, dans \cite{K}, chap. XIV, ou \cite{Tu2}, chap. I), nous obtenons
$$\FT ((\widehat{\D})^{\flat})=\vcenter{\hbox{\input{casgeneral1.pstex_t}}}.$$
Or, si $g\in\Hom_{\R(q)}\left(\FT ( s \otimes k),\FT ( \widetilde{s} )
\right)$, alors
$$g^{\sharp}\in\Hom_{\R(q)}\left(\FT ( s  ),\FT ( \widetilde{s} )
\otimes\FT (k)^*\right)$$
et
$$g^{\sharp}=\vcenter{\hbox{\input{casgeneral2.pstex_t}}}$$
On en d\'eduit que
$$\left(\FT ((\widehat{\D})^{\flat})\right)^{\sharp}=
\vcenter{\hbox{\input{casgeneral3.pstex_t}}}=\vcenter{\hbox{\input{casgeneral4.pstex_t}}}\quad =
(\id_{\FT ( \widetilde{s} )}\otimes\alpha^{\otimes k})\FT 
(\widehat{\D}).$$
Or $(\id_{\FT (\widetilde{s})}\otimes (\alpha^{-1})^{\otimes k}):
\Hom_{\R(q)}\left(\FT (s),\FT (\widetilde{s})\otimes\FT (k)^*\right)
\rightarrow
\Hom_{\R(q)}\left(\FT(s),\FT(t)\right)$ est un isomorphisme, donc 
l'application
$$\Hom_{\V(a)}(s,t)\longrightarrow\Hom_{\R(q)}\left(\FT(s),\FT(t)
\right)$$
qui envoie  $\widehat{\D}$ sur
$\FT (\widehat{\D})=(\id_{\FT(\widetilde{s})}\otimes (\alpha^{-1})^
{\otimes k})\left(\FT((\widehat{\D})^{\flat})\right)^{\sharp}$
est un isomorphisme.\cqfd
\section{Pleine fid\'elit\'e : cas d'une racine de l'unit\'e}
\label{preuvebis}

Nous fixons un nombre complexe $a$ qui une racine primitive $4r^{
\mbox{\scriptsize \`eme}}$ de l'unit\'e avec $r\geq 3$. Posons $q=a^2$.
Rappelons que l'ensemble $J=J(a)$ d\'efini au \S\ref{turaev} est l'ensemble 
$\{ 1,\dots ,r-2\}$. 
Nous allons d\'emontrer que le foncteur $\FB :\VB (a)\rightarrow
\RB (q)$ est pleinement fid\`ele. D'apr\`es le corollaire \ref{fidele}, il 
ne nous reste qu'\`a d\'emontrer la
surjectivit\'e de l'application
\begin{equation}
\label{reste}
\FB:\Hom_{\VB(a)}( s , s ')\longrightarrow\Hom_{\RB(q)}(\FB( s ),\FB(
 s '))
\end{equation}
pour toute paire $(s,s')$ d'objets de $\VB(a)$, ce qui ach\`evera la
d\'emonstration du th\'eor\`eme \ref{principal}.

Remarquons qu'il suffit de v\'erifier que,
pour tout objet $ s $ de $\VB(a)$ et tout entier $k\geq 0$, le 
foncteur $\FB :\Hom_{\VB(a)}( s  ,(k))
\rightarrow\Hom_{\RB(q)}(\FT ( s ),\FT (k))$ 
r\'ealise un isomorphisme d'espaces vectoriels. En effet, la 
d\'emonstration donn\'ee au \S\ref{casgeneral} s'\'etend mot pour mot au
 cas d'une racine de l'unit\'e. Or nous savons que si  
$P$, $Z$ et $W$ sont des $U_q(sl_2)$-modules de dimension finie avec
$Z$ n\'egligeable, alors
$$\hom_{\R(q)}(P\oplus Z,W)\cong \hom_{\R(q)}(P,W),$$ 
(les notations sont celles de \S\ref{mornegl}).
Il suffit donc de consid\'erer les cas o\`u la partie semisimple
 de $\FT ( s  )$ contient un vecteur 
de plus haut poids $k$, 
puisque dans le cas contraire, $\Hom_{\RB(q)}(\FT ( s  ),\FT (k))=0$,
ce qui entra\^ine
$\Hom_{\VB(a)}( s  ,(k))=0$ d'apr\`es le corollaire 
\ref{fidele}.

Nous sommes donc r\'eduits \`a montrer la surjectivit\'e de l'application
\begin{equation}\label{reste2}
\FT:\Hom_{\V(a)}( s ,(k))\longrightarrow\Hom_{\RB(q)}(\FB( s ),
\FB( k))
\end{equation}
pour toute suite $s$ et pour tout plus haut poids $k$ de $\FT(s)$. 
Nous allons proc\'eder par r\'ecurrence sur 
la longueur $m$ de $s$. Si $m=1$, on proc\`ede comme dans le cas 
g\'en\'erique. Au \S\ref{m=2}, nous 
consid\'erons le cas $m=2$.

\subsection{Cas ${\mathbf m=2}$}\label{m=2}

Fixons deux entiers $n_1,n_2\in J$ et $k$ un entier $\geq 0$ tel que 
la partie semisimple de 
$\FT (n_1,n_2)$ contienne un vecteur 
de plus haut poids $k$. D'apr\`es \eqref{CG} et \eqref{CGbis}, ceci 
\'equivaut \`a 
$$k\leq\inf (n_1+n_2,2r-4-n_1-n_2)\quad\text{et}\quad k=n_1+n_2-2j,$$
pour un certain entier $j$ tel que $0\leq j\leq\inf(n_1,n_2)$. 

\begin{prop}
L'application $\FT :\Hom_{\V(a)}((n_1,n_2),(k))\rightarrow\Hom_{\RB(q)}(\FT
 (n_1,n_2),\FT (k))$ est surjective.
\end{prop}
\preuve Au vu du \S\ref{premiercas}, il suffit d'\'etablir que le scalaire
$b_{n_1-i,n_2-i,j-i}\not =0$ pour tout $i=0,\dots ,j-1$.
Remarquons tout d'abord qu'il suffit de le prouver pour $i=0$.
 En effet, le triplet $(n_1-i,n_2-i,j-i)$ v\'erifie les conditions 
impos\'ees \`a $(n_1,n_2,j)$ pour tout $i=0,\dots ,j-1$. Rappelons que 
$$b_{n_1,n_2,j}=q^{-n_2+j-1}\frac{[n_1+n_2-j+1]_q}{[n_1]_q}$$
et que $\inf (n_1+n_2,2r-4-n_1-n_2)=n_1+n_2\Longleftrightarrow n_1+n_2
\leq r-2$.
Nous devons montrer que $[n_1+n_2-j+1]_q\not =0$. Puisque $0<n_1+n_2-j+
1<2r$, il nous reste \`a montrer que 
$n_1+n_2-j+1\not =r$ puisque 
$[n]_q=0$ si et seulement si $n\equiv 0 \,\mbox{mod}\, r$.
Supposons par l'absurde que $n_1+n_2-j+1=r$. Puisque $j\geq 1$, nous 
avons $n_1+n_2>r-2$ et donc
$$\inf (n_1+n_2,2r-4-n_1-n_2)=2r-4-n_1-n_2.$$
Par cons\'equent, $n_1+n_2-2j\leq 2r-4-n_1-n_2$, ce qui est \'equivalent \`a
 $j\geq n_1+n_2-r+2$ et entra\^ine 
une contradiction.\cqfd

\subsection{Cas $\mathbf m\geq 3$}
\label{demo}
On part d'un objet $s$ de longueur $m\geq3$. Soit $s'=(n_1,\dots ,
n_{m-1})$ l'objet de $\V(a)$ d\'efini par
$s=s'\otimes (n_m)$.
Soit $k$ un entier $\geq 0$ tel que la partie semisimple de $\FT ( s  )
$ contienne un vecteur de plus haut poids $k$. 
Par r\'ecurrence, nous admettons que l'application \eqref{reste2} est 
surjective pour tout objet de longueur
$<m$. 

Par analogie avec le cas g\'en\'erique, soit $\{ b_1,\dots ,b_N\}$ l'ensemble des
entiers tel que
pour tout $j=1,\dots ,N$, la dimension $p_j$ de l'espace de plus haut 
poids $b_j$ de la partie semisimple de $\FT(s')$ soit strictement positive et tel 
qu'il existe un 
entier $\ell _j$ v\'erifiant $0\leq \ell _j\leq \inf (b_j,n_m)$ et 
$b_j+n_m-2\ell _j=k$. Comme au \S\ref{deuxiemecas}, les deux lemmes suivant 
justifient la donn\'ee de cet ensemble.

\begin{lem}
\label{nonvide}
Avec les notations pr\'ec\'edente, l'ensemble $\{ b_1,\dots ,b_N\}$ est non vide.
\end{lem} 
\preuve  D\'ecomposons $\FT(s')$ en ses parties semisimples et 
n\'egligeables : $\FT ( s ')\cong P'\oplus Z'$.
Alors
$$\FT ( s )\cong\FT ( s ')\otimes\FT (n_m)\cong (P'\otimes V_{n_m})
\oplus (Z'\otimes V_{n_m}).$$
Puisque $P'\cong \bigoplus_{i=1}^{d}V_{k_i}$ avec $0\leq k_i\leq r-2$
 pour tout $i=1,\dots ,d$,
nous avons
$$P'\otimes V_{n_m}\cong\bigoplus_{i=1}^{d}(V_{k_i}\otimes V_{n_m}).$$
Or, d'apr\`es \eqref{CGbis},
\begin{equation}\label{formule}
V_{k_i}\otimes V_{n_m}\cong\bigoplus_{\stackrel{t=|n_m-k_i|}{
\scriptscriptstyle t\equiv k_i+n_m\,\mbox{\tiny mod}\, 2}}
^{\scriptscriptstyle \inf (k_i+n_m,2r-4-k_i-n_m)}\!\!\!\!\!\!\!
\!\!\!\!\!\!V_t\oplus Z
\end{equation}
o\`u $Z$ est un module n\'egligeable.
De plus, nous savons qu'il  existe un vecteur de plus haut poids 
$k$ dans la partie semisimple de $\FT ( s  )$. 
Nous en d\'eduisons que pour un certain $k_i$, il existe 
un indice $t$ dans \eqref{formule} tel que $t=k$, et donc, d'apr\`es 
\eqref{vecteur+hp}, des entiers $i$
et $\ell $ tels que $k_i+n_m-2\ell =k$ et $0\leq \ell \leq\inf 
(k_i,n_m)$, ce qui prouve que l'ensemble 
$\{b_1,\dots ,b_N\}$ est non vide.\cqfd
Pour tout $j=1,\dots ,N$, soit $w_1^j,\dots ,w_{p_j}^j$ une base de l'espace
de plus haut poids $b_j$ de la partie semisimple de $\FT(s')$. En proc\'edant comme 
pour le lemme \ref{lem12}, on \'etablit le lemme suivant.
\begin{lem}
\label{15}
La famille de vecteurs $v_p^{j}=v(w_p^{j},v_0^{n_m},\ell _j)$
$1\leq j\leq N,\; 1\leq p\leq p_j$, forme une base de l'espace de
plus haut poids $k$ de la partie semisimple
de $\FT(s)$.
\end{lem}

Nous poursuivons comme dans le cas g\'en\'erique.
Toutefois, remarquons que la proposition \ref{egalitedim} n'est pas valable
dans le cas d'une racine de l'unit\'e.

\begin{prop}
\label{16}
L'application 
$\FT :\Hom_{\V (a)}(s,(k))\longrightarrow\Hom_{\RB(q)}(\FB(s),\FB(k))$
est surjective.
\end{prop}
\preuve Etablissons cette proposition par analogie avec la proposition 
\ref{injectivite}, {\em i.e.} montrons qu'il
existe une base $v_1,\dots ,v_{g'}$ de l'espace de plus haut poids 
$k$ de la partie semisimple de $\FT ( s  )$, 
des diagrammes $\D_1,\dots ,\D_{g'}$ du bon type $( s  ,(k))$ et des
 nombres complexes  
$\alpha_1,\dots ,\alpha_{g'},$ tous non nuls tels que l'on ait
\begin{align*}
\F (\D_1)(v_1)&=\alpha_1\, v_0^{\otimes k}, \\
\F (\D_p)(v_{p'})&=0\phantom{_p\, v_0^{\otimes k}}
                       \quad\text{si $1\leq p'<p\leq g'$}, \\
\F (\D_p)(v_p)&=\alpha_p\, v_0^{\otimes k}\quad\text{si 
$p\leq g'$},
\end{align*}
o\`u $g'=\dim\left(\Hom_{\RB(q)}(\FB(s),\FB(k))\right)$.
Puisque la dimension de l'espace de plus haut poids
 $k$ de la partie semisimple de $\FT (s)$
dans le cas d'une racine de l'unit\'e est inf\'erieure \`a la dimension de 
l'espace de plus haut poids de $\FT (s)$ 
dans le cas g\'en\'erique, on a $g'\leq g=\dim\left(\Hom_{\V (a)}(s,(k))\right)$.
La d\'emonstration de la proposition \ref{prop16}
 s'adapte donc au cas pr\'esent d\`es que l'on
remplace les espaces de plus haut poids par les espaces de plus haut
 poids des parties semisimples
 et le lemme \ref{lem12} par le lemme
\ref{15}. Nous conservons la d\'efinition des diagrammes $\overline
{\D_{j,p}}$ et la proposition \ref{recdiag}. \cqfd
A partir de la proposition \ref{16} nous proc\'edons comme dans la
 d\'emonstration du th\'eor\`eme \ref{theoreme}
\`a la fin du \S\ref{deuxiemecas} qui s'applique encore ici pour ce qui 
est de la surjectivit\'e. En conclusion, 
l'application \eqref{reste2} est surjective pour tout $m\geq 3$.

\bigskip\noindent {\bf Remerciements.} Je remercie Christian Kassel 
qui m'a soumis le probl\`eme. 
Son aide a \'et\'e d\'ecisive pour l'\'elaboration de cet article.

\bibliography{sl2}
\bibliographystyle{abbrv}

\end{document}